\documentclass{article}
\usepackage{amssymb}
\usepackage{amsfonts}
\usepackage{amsmath}

\newtheorem{corollary}{Corollary}[section]

\newtheorem{definition}{Definition}[section]
\newtheorem{example}{Example}[section]

\newtheorem{lemma}{Lemma}[section]

\newtheorem{proposition}{Proposition}[section]
\newtheorem{remark}{Remark}[section]

\newtheorem{theorem}{Theorem}[section]
\newenvironment{proof}[1][Proof]{\noindent\textbf{#1.} }{\ \rule{0.5em}{0.5em}}

\input{tcilatex}
\begin{document}

\title{Unital Full Amalgamated Free Products of MF Algebras}
\author{Qihui Li $\ \ \ \ \ \ \ \ \ $Junhao Shen\thanks{%
The research of second author is partially sponsored by an NSF grant.}}
\maketitle

\begin{abstract}
In the paper, we consider the question whether a unital full free product of
MF algebras with amalgamation over a finite dimensional C*-algebra is an MF
algebra. First, we show that, under a natural condition, a unital full free
product of two separable residually finite dimensional (RFD) C*-algebras
with amalgamation over a finite dimensional C*-algebra is again a separable
RFD C*-algebra. Applying this result on MF C*-algebras, we show that, under
a natual condition, a unital full free product of two MF algebras is again
an MF algebra. As an application, we show that a unital full free product of
two AF algebras with amalgamation over an AF algebra is an MF algebra if
there are faithful tracial states on each of these two AF algebras such that
the restricitions on the common subalgebra agree.
\end{abstract}

\section{\protect\bigskip Introduction}

The concept of MF algebras was first introduced by Blackadar and Kirchberg
in \cite{[BK]}. If a separable C*- algebra $\mathcal{A}$ can be embedded
into $\tprod\limits_{k}\mathcal{M}_{n_{k}}\left( \mathbb{C}\right) /\sum_{k}%
\mathcal{M}_{n_{k}}\left( \mathbb{C}\right) $ for a sequence of positive
integers $\left\{ n_{k}\right\} _{k=1}^{\mathcal{1}},$ then $\mathcal{A}$ is
called an MF algebra. Many properties of MF algebras were discussed in \cite%
{[BK]}. For example, it was shown there that an inductive limit of MF
algebras is an MF algebra and every subalgebra of an MF algebra is an MF
algebra. This class of C*-algebras is of interest for many reasons. For
example, it plays an important role in the classification of C$^{\ast }$%
-algebras and it is connected to the question whether the Ext semigroup, in
the sense of Brown, Douglas and Fillmore, of a unital C$^{\ast }$-algebra is
a group (see the striking result of Haagerup and Thorbj{\o }rnsen on $%
Ext(C_{r}^{\ast }(F_{2})$). This notion is also closely connected to
Voiculescu's topological free entropy dimension for a family of self-adjoint
elements $x_{1},\cdots ,x_{n}$ in a unital C*-algebra $\mathcal{A}$ \cite%
{DV5}$.$

The class of MF algebras contains all residually finite-dimensional
C*-algebras and quasidiagonal C*-algebras. Recall that a C*-algebra is said
to be residually finite-dimensional (RFD) if it has a separating family of
finite-dimensional representations. In \cite{CM}, Choi showed that the full
C*-algebra of the free group on two generators is RFD. Later, in [\ref{RT}],
Exel and Loring showed that the unital full free product of two unital RFD
C*-algebras is RFD. In the same paper they gave several equivalent
conditions for the RFD property. In \cite{[ADRL]}, Armstrong, Dykema, Exel
and Li characterized the RFD property of unital full amalgamated free
products of finite dimensional C*-algebras, which extends an earlier result
by Brown and Dykema in [\ref{[BD]}].

Quasidiagonal operators on separable Hilbert spaces were defined by P. R.
Halmos \cite{HA} as compact perturbations of block-diagonal operators. A
generalized notion of quasidiagonal operators to sets of operators is the
concept of quasidiagonal sets of operators. A C*-algebra $\mathcal{A}$ is
quasidiagonal (QD) if there is a faithful representation $\rho $ such that $%
\rho \left( \mathcal{A}\right) $ is a quasidiagonal set of operators. This
class of C*-algebras has been studied for more than 30 years. In [\ref{[BO]}%
], it has been shown that a full free product of two unital QD C*-algebras
amalgamated over units is QD.

In this paper, we consider the question whether a unital full amalgamated
free product of two RFD C*-algebras, or two MF algebras, with amalgamation
over a common C*-algebra is, again, an RFD C*-algebra, or an MF algebra
respectively. One example (see Example 2.1) is given to show that the answer
to this general question is no. Under natural restrictions, we are able to
provide an affirmative answer when we consider a unital full free product of
two MF algebras (or two RFD C*-algebras, two QD C*-algebras) with
amalgamation over a full matrix algebra. Our main results about unital RFD
C*-algebras are as follows:\medskip 

\textbf{Theorem 4.1 }Let $\mathcal{A}$, $\mathcal{B}$ be unital RFD algebras
and $\mathcal{D}$ be a finite-dimensional C*-algebra. Let $\psi _{\mathcal{A}%
}:\mathcal{D\rightarrow A}$ and $\psi _{\mathcal{B}}:$ $\mathcal{%
D\rightarrow B}$ be unital embeddings. Let $\mathcal{A}\underset{\mathcal{D}}%
{\mathcal{\ast }}\mathcal{B}$ be the corresponding unital full amalgamated
free product. Then $\mathcal{A}\underset{\mathcal{D}}{\mathcal{\ast }}%
\mathcal{B}$ is RFD if and only if there is a sequence $\left\{
k_{n}\right\} _{n=1}^{\mathcal{1}}$ of positive integers with unital
embeddings $q_{1}:\mathcal{A\rightarrow }\prod_{n=1}^{\infty }\mathcal{M}%
_{k_{n}}(\mathbb{C})$ and $q_{2}:\mathcal{B\rightarrow }\prod_{n=1}^{\infty }%
\mathcal{M}_{k_{n}}(\mathbb{C})$ such that the following diagram commutes

\begin{equation*}
\begin{array}{ccc}
\mathcal{D} & \overset{\psi _{\mathcal{A}}}{\hookrightarrow } & \mathcal{A}
\\ 
^{\psi _{\mathcal{B}}}\downarrow &  & {\Huge \downarrow }^{q_{1}} \\ 
\mathcal{B} & \overset{q_{2}}{\hookrightarrow } & \prod_{m=1}^{\infty }%
\mathcal{M}_{k_{m}}(\mathbb{C)}%
\end{array}%
\end{equation*}%
\medskip

\textbf{Corollary 4.1 }Suppose that $\mathcal{A}$ is RFD and $\mathcal{D}$
is a unital finite-dimensional C*-subalgebra of $\mathcal{A}.$ Then $%
\mathcal{A\ast }_{\mathcal{D}}\mathcal{A}$ is RFD.\medskip

\bigskip \textbf{Theorem 5.1 }Let $\mathcal{A}$ and $\mathcal{B}$ be unital
MF-algebras and $\mathcal{D}$ be a finite-dimensional C*-algebra. Let $\psi
_{1}:\mathcal{D\rightarrow A}$ and $\psi _{2}:$ $\mathcal{D\rightarrow B}$
be unital embeddings. Then $\mathcal{A}\underset{\mathcal{D}}{\mathcal{\ast }%
}\mathcal{B}$ is an MF algebra if and only if there is a sequence $\left\{
k_{n}\right\} _{n=1}^{\mathcal{1}}$ of integers with unital embeddings $%
q_{1}:\mathcal{A\rightarrow }\prod_{n=1}^{\infty }\mathcal{M}_{k_{n}}(%
\mathbb{C})/\overset{\infty }{\underset{n=1}{\sum }}\mathcal{M}_{k_{n}}(%
\mathbb{C})$ and $q_{2}:\mathcal{B\rightarrow }\prod_{n=1}^{\infty }\mathcal{%
M}_{k_{n}}(\mathbb{C})/\overset{\infty }{\underset{n=1}{\sum }}\mathcal{M}%
_{k_{n}}(\mathbb{C})$ such that the following diagram commutes

\begin{equation*}
\begin{array}{ccc}
\mathcal{D} & \overset{\psi _{\mathcal{A}}}{\hookrightarrow } & \mathcal{A}
\\ 
^{\psi _{\mathcal{B}}}\downarrow &  & {\Huge \downarrow }^{q_{1}} \\ 
\mathcal{B} & \overset{q_{2}}{\hookrightarrow } & \prod_{m=1}^{\infty }%
\mathcal{M}_{k_{m}}(\mathbb{C)}/\tsum \mathcal{M}_{k_{m}}(\mathbb{C)}%
\end{array}%
\end{equation*}

\bigskip

\bigskip

When the common C*-algebra $\mathcal{D}$ is an AF algebra, we obtain
following results. \medskip

\textbf{Theorem 5.3 }Suppose that $\mathcal{A}\supset \mathcal{D}\subset 
\mathcal{B}$ are unital inclusions of unital MF algebras where $\mathcal{D}$
is an AF algebra. Then the unital full free product $\mathcal{A}\ast _{%
\mathcal{D}}\mathcal{B}$ of $\mathcal{A}$ and $\mathcal{B}$ with
amalgamation over $\mathcal{D}$ is an MF algebra if and only if there is a
sequence of positive integers $\{k_{n}\}_{n=1}^{\infty }$ such that the
following diagram is commutative 
\begin{equation*}
\begin{array}{ccc}
\mathcal{D} & \subseteq & \mathcal{A} \\ 
\cap &  & \cap \\ 
\mathcal{B} & \subseteq & \tprod\limits_{n}\mathcal{M}_{k_{n}}\left( \mathbb{%
C}\right) /\sum_{n}\mathcal{M}_{k_{n}}\left( \mathbb{C}\right)%
\end{array}%
\end{equation*}

\textbf{Theorem 5.4 \ }Suppose that $\mathcal{A}\supset \mathcal{D}\subset 
\mathcal{B}$ are unital inclusions of AF C$^{\ast }$-algebras. If there are
faithful tracial states $\tau _{\mathcal{A}}$ and $\tau _{\mathcal{B}}$ on $%
\mathcal{A}$ and $\mathcal{B}$ respectively, such that 
\begin{equation*}
\tau _{\mathcal{A}}(x)=\tau _{\mathcal{B}}(x),\qquad \forall \ x\in \mathcal{%
D},
\end{equation*}%
then $\mathcal{A}\ast _{\mathcal{D}}\mathcal{B}$ is an MF algebra.\medskip

\textbf{Corollary 5.2 \ }Suppose that $\mathcal{A\supseteq D\subseteq B}$
are unital inclusions of C*-algebras where $\mathcal{A}$, $\mathcal{B}$ are
UHF algebras and $\mathcal{D}$ is an AF algebra. Then $\mathcal{A}\underset{%
\mathcal{D}}{\mathcal{\ast }}\mathcal{B}$ is an MF algebra if and only if 
\begin{equation*}
\tau _{\mathcal{A}}\left( z\right) =\tau _{\mathcal{B}}\left( x\right) \ 
\text{for\ each\ }z\in \mathcal{D},
\end{equation*}%
where $\tau _{\mathcal{A}}$ and $\tau _{\mathcal{B}}$ are faithful tracial
states on $\mathcal{A}$ and $\mathcal{B}$ respectively.

A brief overview of this paper is as follows. In Section 2, we recall the
definitions of MF algebras and unital full amalgamated free product of
unital C*-algebas. One example is given at the end of the section to show
that a unital full amalgamated free product of unital MF (or RFD, QD)
algebras may not be MF (or RFD, QD) again. Section 3 contains definitions
and some basic properties of quasidiagonal operators, quasidiagonal sets of
operators and QD C*-algebras. This section is designed to make the paper
self-contained. Section 4 is devoted to results on the full amalgamated free
products of unital RFD C*-algebras. In Section 5, we first consider unital
full free products of unital MF algebras with amalgamation over
finite-dimensional C*-subalgebras. Then we consider the case when a common
unital C*-subalgebra is an infinite-dimensional C*-algebra. More precisely,
we consider a case when the common unital C*-subalgebra is an AF algebra.

\section{Definitions and Preliminaries}

\subsection{Noncommutative polynomials}

In this article, we always assume that all C*-algebras are unital separable C%
$^{\ast }$-algebras. We use notation C*$\left( x_{1},x_{2,}\cdots \right) $
to denote the unital C*-algebra generated by $\left\{ x_{1},x_{2},\cdots
\right\} .$ Let $\mathbb{C\langle }\mathbf{X}_{1},\ldots ,\mathbf{X}_{n}%
\mathbb{\rangle }$ be the set of all noncommutative polynomials in the
indeterminants $\mathbf{X}_{1},\ldots ,\mathbf{X}_{n}$. Let $\mathbb{C}_{%
\mathbb{Q}}=\mathbb{Q}+i\mathbb{Q}$ denote the complex-rational numbers,
i.e., the numbers whose real and imaginary parts are rational. Then the set $%
\mathbb{C}_{\mathbb{Q}}\mathbb{\langle }\mathbf{X}_{1},\ldots ,\mathbf{X}_{n}%
\mathbb{\rangle }$ of noncommutative polynomials with complex-rational
coefficients is countable. Throughout this paper we write

\begin{equation*}
\mathbb{C\langle}\mathbf{X}_{1},\mathbf{X}_{2},\cdots\mathbb{\rangle=\cup }%
_{m=1}^{\infty}\mathbb{C\langle}\mathbf{X}_{1},\mathbf{X}_{2},\cdots \mathbf{%
X}_{m}\mathbb{\rangle}
\end{equation*}
and 
\begin{equation*}
\mathbb{C}_{\mathbb{Q}}\mathbb{\langle}\mathbf{X}_{1},\mathbf{X}_{2},\cdots%
\mathbb{\rangle=\cup}_{m=1}^{\infty}\mathbb{C}_{\mathbb{Q}}\mathbb{\langle}%
\mathbf{X}_{1},\mathbf{X}_{2},\cdots\mathbf{X}_{m}\mathbb{\rangle}.
\end{equation*}

Let $\left\{ P_{r}\right\} _{r=1}^{\mathcal{1}}$ be the collection of all
noncommutative polynomials in $\mathbb{C}_{\mathbb{Q}}\mathbb{\langle }%
\mathbf{X}_{1},\mathbf{X}_{2},\cdots \mathbb{\rangle }$ with rational
complex coefficients.

\subsection{Blackadar and Kirchberg's MF Algebras}

We will recall an equivalent definition of MF algebras and some basic
properties. Let us fix notations first. We assume that $\mathcal{H}$ is a
separable complex Hilbert space and $\mathcal{B}(\mathcal{H})$ is the set of
all bounded operators on $\mathcal{H}$. Suppose $\{x,x_{k}\}_{k=1}^{\infty }$
is a family of elements in $\mathcal{B}(\mathcal{H})$. We say $%
x_{k}\rightarrow x$ in $\ast $-SOT ($\ast $-strong operator topology) if and
only if $x_{k}\rightarrow x$ in SOT and $x_{k}^{\ast }\rightarrow x^{\ast }$
in SOT. Suppose $\{x_{1},\ldots ,x_{n}\}$ and $\{x_{1}^{(k)},\ldots
,x_{n}^{(k)}\}_{k=1}^{\infty }$ are families of elements in $\mathcal{B}(%
\mathcal{H})$. We say 
\begin{equation*}
\langle x_{1}^{(k)},\ldots ,x_{n}^{(k)}\rangle \rightarrow \langle
x_{1},\ldots ,x_{n}\rangle \ in\ \ast \text{-}SOT,\ \ as\ k\rightarrow \infty
\end{equation*}%
if and only if 
\begin{equation*}
x_{i}^{(k)}\rightarrow x_{i}\ in\ \ast \text{-}SOT,\ \ as\ k\rightarrow
\infty ,\ \ \forall \ 1\leq i\leq n.
\end{equation*}%
Suppose $\mathcal{A}$ is a separable unital C*-algebra on a Hilbert space $%
\mathcal{H}$. Let $\mathcal{H}^{\infty }=\oplus _{\mathbb{N}}\mathcal{H}$
and, for any $x\in \mathcal{A}$, let $x^{\infty }$ be the element $\oplus _{%
\mathbb{N}}x=(x,x,x,\ldots )$ in $\prod_{k\in \mathbb{N}}\mathcal{A}%
^{(k)}\subset B(\mathcal{H}^{\infty })$, where $\mathcal{A}^{(k)}$ is the $k$%
-th copy of $\mathcal{A}$.

\bigskip Suppose $\{\mathcal{M}_{k_{n}}(\mathbb{C})\}_{n=1}^{\infty }$ is a
sequence of complex matrix algebras. We introduce the C*-direct product $%
\prod_{m=1}^{\infty }\mathcal{M}_{k_{m}}(\mathbb{C)}$ of $\{\mathcal{M}%
_{k_{n}}(\mathbb{C})\}_{n=1}^{\infty }$ as follows: 
\begin{equation*}
\prod_{n=1}^{\infty }\mathcal{M}_{k_{n}}(\mathbb{C})=\{(Y_{n})_{n=1}^{\infty
}\ |\ \forall \ n\geq 1,\ Y_{n}\in \mathcal{M}_{k_{n}}(C)\ \text{ and }\
\left\Vert (Y_{n})_{n=1}^{\infty }\right\Vert =\sup_{n\geq 1}\Vert
Y_{n}\Vert <\infty \}.
\end{equation*}%
Furthermore, we can introduce a norm-closed two sided ideal in $%
\prod_{n=1}^{\infty }\mathcal{M}_{k_{n}}(\mathbb{C})$ as follows: 
\begin{equation*}
\overset{\infty }{\underset{n=1}{\sum }}\mathcal{M}_{k_{n}}(\mathbb{C}%
)=\left\{ \left( Y_{n}\right) _{n=1}^{\infty }\in \prod_{n=1}^{\infty }%
\mathcal{M}_{k_{n}}(\mathbb{C}):\lim\limits_{n\rightarrow \infty }\left\Vert
Y_{n}\right\Vert =0\right\} .
\end{equation*}%
Let $\pi $ be the quotient map from $\prod_{n=1}^{\infty }\mathcal{M}%
_{k_{n}}(\mathbb{C})$ to $\prod_{n=1}^{\infty }\mathcal{M}_{k_{n}}(\mathbb{C}%
)/\overset{\infty }{\underset{n=1}{\sum }}\mathcal{M}_{k_{n}}(\mathbb{C})$.
Then 
\begin{equation*}
\prod_{n=1}^{\infty }\mathcal{M}_{k_{n}}(\mathbb{C})/\overset{\infty }{%
\underset{n=1}{\sum }}\mathcal{M}_{k_{n}}(\mathbb{C})
\end{equation*}%
is a unital C*-algebra. If we denote $\pi \left( \left( Y_{n}\right)
_{n=1}^{\infty }\right) $ by $\left[ \left( Y_{n}\right) _{n}\right] $, then 
\begin{equation*}
\left\Vert \left[ \left( Y_{n}\right) _{n}\right] \right\Vert =\underset{%
n\rightarrow \mathcal{1}}{\lim \sup }\left\Vert Y_{n}\right\Vert .
\end{equation*}%
Now we are ready to recall an equivalent definition of MF algebras given by
Blackadar and Kirchberg [\ref{[BK]}].

\begin{definition}
\label{MF}(Theorem 3.2.2, [\ref{[BK]}]) Let $\mathcal{A}$ be a separable
C*-algebra. If $\mathcal{A}$ can be embedded as a C*-subalgebra of $%
\prod_{n=1}^{\infty }\mathcal{M}_{k_{n}}(\mathbb{C})/\overset{\infty }{%
\underset{n=1}{\sum }}\mathcal{M}_{k_{n}}(\mathbb{C})$ for a sequence $%
\left\{ k_{n}\right\} _{n=1}^{\mathcal{1}}$ of integers$,$ then $\mathcal{A}$
is called an MF algebra.
\end{definition}

The following Theorem is one of the key ingredients for showing our main
results in this paper.

\begin{theorem}
\label{100}(Theorem 5.1.2, \cite{HLS}) Suppose that $\mathcal{A}$ is a
unital C$^{\ast }$-algebra generated by a sequence of self-adjoint elements $%
x_{1},x_{2},\cdots $ in $\mathcal{A}.$ Then the following are equivalent:

\begin{enumerate}
\item $\mathcal{A}$ is an MF algebra

\item For each $n\in \mathbb{N},$ there are a sequence of positive integers $%
\{m_{k}\}_{k=1}^{\infty }$ and self-adjoint matrices $A_{1}^{(k)},\ldots
,A_{n}^{(k)}$ in $\mathcal{M}_{m_{k}}^{s.a.}(\mathbb{C)}$ for $k=1,2,\ldots $%
, such that, $\forall \ P\in \mathbb{C}\langle X_{1},\ldots ,X_{n}\rangle $, 
\begin{equation*}
\lim_{k\rightarrow \infty }\Vert P(A_{1}^{(k)},\ldots ,A_{n}^{(k)})\Vert
=\Vert P(x_{1},\ldots ,x_{n})\Vert ,
\end{equation*}%
where $\mathbb{C}\langle X_{1},\ldots ,X_{n}\rangle $ is the set of all
noncommutative polynomials in the indeterminants $X_{1},\ldots ,X_{n}$.
\end{enumerate}
\end{theorem}

The examples of MF algebras contain all finite dimensional C*-algebras, AF
(approximately finite dimensional) algebras and quasidiagonal C*-algebras.
In [\ref{HT}], Haagerup and Thorbj$\phi $rnsen showed that $C_{r}^{\ast
}\left( F_{n}\right) $ is an MF algebra for $n\geq 2.$ For more examples of
MF algebras, we refer the reader to [\ref{[BK]}] and [\ref{DS}].

\subsection{Definition of Full Amalgamated Free Products of Unital
C*-algebras}

Recall the definition of full amalgamated free product of unital C*-algebras
as follows.

\begin{definition}
Given C*-algebras $\mathcal{A}$, $\mathcal{B}$ and $\mathcal{D}$ with unital
embeddings (injective $\ast $ -homomorphisms) $\psi _{\mathcal{A}}:\mathcal{%
D\rightarrow A}$ and $\psi _{\mathcal{B}}:\mathcal{D\rightarrow B}$, the
corresponding full amalgamated free product C*-algebra is the C*-algebra $%
\mathcal{C}$, equipped with unital embeddings $\sigma _{\mathcal{A}}:%
\mathcal{A\rightarrow C}$ and $\sigma _{\mathcal{B}}:\mathcal{B\rightarrow C}
$ such that $\sigma _{\mathcal{A}}\circ \psi _{\mathcal{A}}=\sigma _{%
\mathcal{B}}\circ \psi _{\mathcal{B}},$ such that $\mathcal{C}$ is generated
by $\sigma _{\mathcal{A}}\left( \mathcal{A}\right) \cup \sigma _{\mathcal{B}%
}\left( \mathcal{B}\right) $ and satisfying the universal property that
whenever $\mathcal{E}$ is a C*-algebra and $\pi _{\mathcal{A}}:\mathcal{%
A\rightarrow E}$ and $\pi _{\mathcal{B}}:\mathcal{B\rightarrow E}$ are $\ast 
$-homomorphisms satisfying $\pi _{\mathcal{A}}\circ \psi _{\mathcal{A}}=\pi
_{\mathcal{B}}\circ \psi _{\mathcal{B}},$ there is a $\ast $-homomorphism $%
\pi :\mathcal{C\rightarrow E}$ such that $\pi \circ \sigma _{\mathcal{A}%
}=\pi _{\mathcal{A}}$ and $\pi \circ \sigma _{\mathcal{B}}=\pi _{\mathcal{B}%
}.$ The full amalgamated free product C*-algebra $\mathcal{C}$ is commonly
denoted by $\mathcal{A}\underset{\mathcal{D}}{\mathcal{\ast }}\mathcal{B}.$
\end{definition}

When $D=\mathbb{C}I,$\ the above definition is the unital full free product $%
\mathcal{A}\ast _{\mathbb{C}}\mathcal{B}$ of $\mathcal{A}$ and $\mathcal{B}$%
{\large . }The following theorem is certainly well-known, and we will give
its proof for the purpose of completeness.

\begin{theorem}
\label{32}Suppose that $\mathcal{A},$ $\mathcal{B}$ and $\mathcal{D}$ are
unital C*-algebras. Then 
\begin{equation*}
\left( \mathcal{A\otimes }_{\max }\mathcal{D}\right) \underset{\mathcal{D}}{%
\mathcal{\ast }}\left( \mathcal{B\otimes }_{\max }\mathcal{D}\right) 
\mathcal{\cong }\left( \mathcal{A}\underset{\mathbb{C}}{\mathcal{\ast }}%
\mathcal{B}\right) \otimes _{\max }\mathcal{D}.
\end{equation*}
\end{theorem}

\begin{proof}
From the definition of unital full free product, we can get two natural
unital embeddings 
\begin{equation*}
\pi _{1}:\mathcal{A\otimes }_{\max }\mathcal{D\rightarrow }\left( \mathcal{A}%
\underset{\mathbb{C}}{\mathcal{\ast }}\mathcal{B}\right) \otimes _{\max }%
\mathcal{D}
\end{equation*}%
and 
\begin{equation*}
\pi _{2}:\mathcal{B\otimes }_{\max }\mathcal{D\rightarrow }\left( \mathcal{A}%
\underset{\mathbb{C}}{\mathcal{\ast }}\mathcal{B}\right) \otimes _{\max }%
\mathcal{D}
\end{equation*}%
from $\mathcal{A\otimes }_{\max }\mathcal{D}$ and $\mathcal{B\otimes }_{\max
}\mathcal{D}$ into $\left( \mathcal{A}\underset{\mathbb{C}}{\mathcal{\ast }}%
\mathcal{B}\right) \otimes _{\max }\mathcal{D}$ respectively. It is clear
that the restrictions of $\pi _{1}$ and $\pi _{2}$ on $I\otimes \mathcal{D}$
agree, i.e., $\pi _{1}|_{I\otimes \mathcal{D}}=\pi _{2}|_{I\otimes \mathcal{D%
}}.$ Suppose $\mathcal{K}$ is a C*-algebra acting on a Hilbert space $%
\mathcal{H}$ such that there are two *-homomorphisms $q_{1}:\mathcal{%
A\otimes }_{\max }\mathcal{D\rightarrow K}$ and $q_{2}:\mathcal{B\otimes }%
_{\max }\mathcal{D\rightarrow K}$ satisfying $q_{1}|_{I\otimes \mathcal{D}%
}=q_{2}|_{I\otimes \mathcal{D}}.$ It implies that $q_{1}\left( \mathcal{%
A\otimes I}\right) $ commutes with $q_{1}\left( I\otimes \mathcal{D}\right) $
in $\mathcal{K}$ and $q_{2}\left( \mathcal{B\otimes I}\right) $ commutes
with $q_{2}\left( I\otimes \mathcal{D}\right) $ in $\mathcal{K}.$ Let $%
\mathcal{M=K\cap }\left( q_{1}\left( 1\otimes \mathcal{D}\right) \right)
^{\prime }=\mathcal{K\cap }\left( q_{2}\left( 1\otimes \mathcal{D}\right)
\right) ^{\prime }.$ Since $q_{1}\left( \mathcal{A\otimes I}\right) $ and $%
q_{2}\left( \mathcal{B\otimes I}\right) $ are subalgebras of $\mathcal{M},$
there is a *-homomorphism $\widetilde{q}:\mathcal{A}\underset{\mathbb{C}}{%
\mathcal{\ast }}\mathcal{B\rightarrow M}$ by the definition of unital full
free product. Moreover, the image $\widetilde{q}\left( \mathcal{A\ast }_{%
\mathbb{C}}\mathcal{B}\right) $ of $\mathcal{A\ast }_{\mathbb{C}}\mathcal{B}$
under $\widetilde{q}$ commutes with $q_{1}\left( I\otimes \mathcal{D}\right) 
$ in $\mathcal{K}.$ From the definition of maximal C*-norm on tensor product
of two C*-algebras, there is a *-homomorphism 
\begin{equation*}
q:\left( \mathcal{A}\underset{\mathbb{C}}{\mathcal{\ast }}\mathcal{B}\right)
\otimes _{\max }\mathcal{D\rightarrow K}.
\end{equation*}%
such that $q\circ \pi _{1}=q_{1}$ and $q\circ \pi _{2}=q_{2}.$ The desired
conclusion now follows from the definition of full amalgamated free products
of unital C*-algebras.
\end{proof}

Combining the following Theorem and preceding result, we are able to obtain
our first result about unital full amalgamated free product of MF algebras,
that is \ Proposition \ref{5}.

\begin{theorem}
(Theorem 5.1.4., [\ref{[HLS]}])\label{13.4} Suppose $\left\{ A_{i}:i\in 
\mathbb{N}\right\} $ is a countable family of separable MF C*-algebras. Then
the unital full free product $\mathcal{A=\ast }_{\mathbb{C}}\mathcal{A}_{i}$
is an MF algebra.
\end{theorem}

\begin{proposition}
\label{5}Let $\mathcal{A}$ and $\mathcal{B}$ be separable unital
C*-algebras. If $\mathcal{D}$ can be embedded as a unital C*-subalgebra of $%
\mathcal{A}$ and $\mathcal{B}$ respectively, and $\mathcal{D}$ is
*-isomorphic to a full matrix algebra $\mathcal{M}_{n}\left( \mathbb{C}%
\right) $ for some integer $n,$ then the full amalgamated free product $%
\mathcal{A}\underset{\mathcal{D}}{\mathcal{\ast }}\mathcal{B}$ is an MF
algebra if and only if $\mathcal{A}$ and $\mathcal{B}$ are both MF.
\end{proposition}

\begin{proof}
If $\mathcal{A}\underset{\mathcal{D}}{\mathcal{\ast }}\mathcal{B}$ is a
unital MF algebra, then it is easy to see that $\mathcal{A}$ and $\mathcal{B}
$ are MF. On the other hand, since $\mathcal{D}$ is *-isomorphic to a full
matrix algebra$,$ from Lemma 6.6.3 in \cite{KR}, it follows that $\mathcal{%
A\cong A}^{\prime }\otimes \mathcal{D}$ and $\mathcal{B\cong B}^{\prime
}\otimes \mathcal{D}$ where $\mathcal{A}^{\prime }$ and $\mathcal{B}^{\prime
}$ are C*-subalgebras of $\mathcal{A}$ and $\mathcal{B}$ respectively.
Therefore $\mathcal{A}^{\prime }$ and $\mathcal{B}^{\prime }$ are MF as
well. Then the desired conclusion follows from Theorem \ref{13.4} and \ref%
{32}.
\end{proof}

Recall that a separable C*-algebra $\mathcal{R}$ is said to be residually
finite-dimensional (RFD) if for each $x$ $\in \mathcal{R}$ there exists a
*-homomorphism $\pi :\mathcal{R\rightarrow B}$ such that $\dim \left( 
\mathcal{B}\right) <\infty $ and $\pi \left( x\right) \neq 0.$ It is easy to
see that such algebras have a faithful representation whose image is a block
diagonal set of operators. The Proposition \ref{5} is stated for unital MF
algebras, but same conclusion holds when we consider unital RFD C*-algebras
or unital quasidiagonal C*-algebras (we will recall the definition of
quasidiagonal C*-algebra later).

\begin{lemma}
\label{13}(Theorem 4.2.,[\ref{[ADRL]}]) Consider unital inclusions of
C*-algebras $\mathcal{A\supseteq D\subseteq B}$ with $\mathcal{A}$ and $%
\mathcal{B}$ finite dimensional. Let $\mathcal{A}\underset{\mathcal{D}}{%
\mathcal{\ast }}\mathcal{B}$ be the corresponding full amalgamated free
product. Then $\mathcal{A}\underset{\mathcal{D}}{\mathcal{\ast }}\mathcal{B}$
is residually finite dimensional if and only if there are faithful tracial
states $\tau _{\mathcal{A}}$ on $\mathcal{A}$ and $\tau _{\mathcal{B}}$ on $%
\mathcal{B}$ whose restrictions on $\mathcal{D}$ agree.
\end{lemma}

\begin{remark}
Combining Lemma \ref{13} and the fact that each RFD C*-algebra has a
faithful tracial state, it is not hard to see that $\mathcal{A}\underset{%
\mathcal{D}}{\mathcal{\ast }}\mathcal{B}$ is RFD if and only if $\mathcal{A}%
\underset{\mathcal{D}}{\mathcal{\ast }}\mathcal{B}$ has a faithful tracial
state.
\end{remark}

In [\ref{RT}], Exel and Loring showed that the unital full free product of
two RFD C*-algebra is RFD, which extends an earlier result by Choi in [\ref%
{CM}]. In [\ref{[BO]}], Boca showed that the unital full free product of two
quasidiagonal C*-algebras is also quasidiagonal. The analogous result for MF
algebras is given in [\ref{[HLS]}]. Lemma \ref{13} characterizes when a full
amalgamated free product of finite dimensional C*-algebras is RFD, which
extends an earlier result by Brown and Dykema in [\ref{[BD]}]. So from
Proposition \ref{5} and the above results, it is natural to ask whether a
full amalgamated free product of unital MF (or RFD, quasidiagonal) algebras
is always MF (or RFD, quasidiagonal). For the case when $\mathcal{D}$ is
*-isomorphic to a full matrix algebra, we know that $\mathcal{A}\underset{%
\mathcal{D}}{\mathcal{\ast }}\mathcal{B}$ is MF (or RFD, quasidiagonal) if
and only if $\mathcal{A}$ and $\mathcal{B}$ are both MF (or RFD,
quasidiagonal) by Proposition \ref{5}. But the following example shows that
a full amalgamated free product of two MF (or RFD, quasidiagonal) algebras
may not be MF (or RFD, quasidiagonal) again, even for two full matrix
algebras with amalgamation over the two dimensional C*-algebra $\mathbb{%
C\oplus C}$.

\begin{example}
Let $\mathcal{D=}\mathbb{C\oplus C}$. Suppose $\varphi _{1}:\mathcal{%
D\rightarrow M}_{2}\left( \mathbb{C}\right) $ and $\varphi _{2}:\mathcal{%
D\rightarrow M}_{3}\left( \mathbb{C}\right) $ are unital embeddings such
that 
\begin{equation*}
\varphi _{1}\left( 1\oplus 0\right) =\left( 
\begin{array}{cc}
1 & 0 \\ 
0 & 0%
\end{array}%
\right) \text{ and }\varphi _{2}\left( 1\oplus 0\right) =\left( 
\begin{array}{ccc}
1 & 0 & 0 \\ 
0 & 0 & 0 \\ 
0 & 0 & 0%
\end{array}%
\right)
\end{equation*}%
Then $\mathcal{M}_{2}\left( \mathbb{C}\right) \underset{\mathcal{D}}{%
\mathcal{\ast }}\mathcal{M}_{3}\left( \mathbb{C}\right) $ is not an MF
algebra (therefore it is not RFD or quasidiagonal). Actually, if we assume
that $\mathcal{M}_{2}\left( \mathbb{C}\right) \underset{\mathcal{D}}{%
\mathcal{\ast }}\mathcal{M}_{3}\left( \mathbb{C}\right) $ is an MF algebra,
then there exists a tracial state $\tau $ on $\mathcal{M}_{2}\left( \mathbb{C%
}\right) \underset{\mathcal{D}}{\mathcal{\ast }}\mathcal{M}_{3}\left( 
\mathbb{C}\right) .$ So the restrictions of $\tau $ on $\mathcal{M}%
_{2}\left( \mathbb{C}\right) $ and $\mathcal{M}_{3}\left( \mathbb{C}\right) $
are the unique tracial states on $\mathcal{M}_{2}\left( \mathbb{C}\right) $
and $\mathcal{M}_{3}\left( \mathbb{C}\right) $ respectively. It follows that 
$\tau \left( \varphi _{1}\left( 1\oplus 0\right) \right) =\frac{1}{2}\neq $ $%
\tau $ $\left( \varphi _{2}\left( 1\oplus 0\right) \right) =\frac{1}{3}$
which contradicts to the fact that $\varphi _{1}\left( 1\oplus 0\right)
=\varphi _{2}\left( 0\oplus 1\right) $ in $\mathcal{M}_{2}\left( \mathbb{C}%
\right) \underset{\mathcal{D}}{\mathcal{\ast }}\mathcal{M}_{3}\left( \mathbb{%
C}\right) $. Therefore $\mathcal{M}_{2}\left( \mathbb{C}\right) \underset{%
\mathcal{D}}{\mathcal{\ast }}\mathcal{M}_{3}\left( \mathbb{C}\right) $ is
not MF.
\end{example}

\section{Basic Properties of Quasidiagonal Algebras}

The proof of one of our main theorems is based on the understanding of
quasidiagonal C*-algebras. Therefore, we will recall some results about
quasidiagonal C*-algebras for the reader's convenience, but we will not give
their proofs. We refer the reader to [\ref{[B]}] for a comprehensive
treatment of this important class of C*-algebras.

\begin{definition}
\label{6}A subset $\Omega\subset\mathcal{B}\left( \mathcal{H}\right) $ is
called a \emph{quasidiagonal set of operators} if for each finite set $%
\mathcal{S}\subset\Omega,$ finite set $F\subset\mathcal{H}$ and $%
\varepsilon>0$ there exists a finite rank projection $P\in\mathcal{B(H)}$
such that $\left\Vert SP-PS\right\Vert \leq\varepsilon$ and $\left\Vert
Px-x\right\Vert \leq\varepsilon$ for all $S\in\mathcal{S}$ and $x\in F.$
\end{definition}

\begin{definition}
\label{7}A C*-algebra $\mathcal{A}$ is called \emph{quasidiagonal (QD)} if
there exists a faithful representation $\pi :\mathcal{A\rightarrow B(H)}$
such that $\pi \left( \mathcal{A}\right) $ is a quasidiagonal set of
operators.
\end{definition}

\begin{definition}
\label{8}Let $\pi:$ $\mathcal{A}\rightarrow\mathcal{B(H)}$ be a faithful
representation of a C*-algebra $\mathcal{A}.$ Then $\pi$ is called \emph{%
essential} if $\pi(\mathcal{A)}$ contains no nonzero finite rank operators.
\end{definition}

The next lemma is a fundamental result about representations of
quasidiagonal C*-algebras.

\begin{lemma}
(Theorem 1.7, [\ref{Vo}])\label{9} Let $\pi :\mathcal{A}\rightarrow \mathcal{%
B(H)}$ be a faithful essential representation. Then $\mathcal{A}$ is
quasidiagonal if and only if $\pi (\mathcal{A)}$ is a quasidiagonal set of
operators.
\end{lemma}

The following lemma is an important ingredient in the proof of Proposition %
\ref{15}.

\begin{lemma}
\label{10}(Lemma 2.1, \cite{DS}) Suppose that $\mathcal{A\subset B(H)}$ is a
separable unital quasidiagonal C*-algebra and $x_{1},\cdots ,x_{n}$ are
self-adjoint elements in$\mathcal{\ A}.$ For any $\varepsilon >0,$ any
finite subset $\left\{ f_{1},\cdots ,f_{r}\right\} $ of $\mathbb{C}%
\left\langle \mathbf{X}_{1},\cdots ,\mathbf{X}_{n}\right\rangle $ and any
finite subset $\left\{ \xi _{1},\cdots ,\xi _{r}\right\} $ of $\mathcal{H}$,
there is a finite rank projection $p$ in $\mathcal{B(H)}$ such that:

(i) $\left\Vert \xi _{k}-p\xi _{k}\right\Vert <\varepsilon ,$ $\left\Vert
\left( px_{i}p-x_{i}\right) \xi _{k}\right\Vert <\varepsilon ,$ for all $%
1\leq i\leq n$ and $1\leq k\leq r;$

(ii) $\left\vert \left\Vert f_{j}\left( px_{1}p,\cdots ,px_{n}p\right)
\right\Vert _{\mathcal{B(}p\mathcal{H)}}-\left\Vert f_{j}\left( x_{1},\cdots
,x_{n}\right) \right\Vert \right\vert <\varepsilon ,$ for all $1\leq j\leq
r. $
\end{lemma}

Using Lemma \ref{10}, it is easy to see that all quasidiagonal C*-algebras
are MF algebras.

\begin{lemma}
\label{11}(Proposition 7.4, \cite[7.4]{[B]}) If $\left\{ A_{n}\right\} $ is
a sequence of C*-algebras then $\tprod\limits_{n\in \mathbb{N}}A_{n}$ is QD
if and only if each $A_{n}$ is QD.
\end{lemma}

The examples of quasidiagonal C*-algebras include all abelian C*-algebras
and finite-dimensional C*-algebras as well as residually finite-dimensional
C*-algebras.

\section{Full Amalgamated Free Products of RFD C*-algebras}

First, we will give the following well-known lemma. For completeness, we
include the proof.

\begin{lemma}
\label{finite projections}Given $0<\epsilon <1$ and $n\in \mathbb{N}.$ For
any two families of $n$ pairwise orthogonal projections $\left\{
P_{1},\cdots ,P_{n}\right\} $ and $\left\{ Q_{1},\cdots ,Q_{n}\right\} $ in $%
n$dimensional unital abelian C*-subalgebras $\mathcal{A}$ and $\mathcal{B}$
in $\mathcal{B}\left( \mathcal{H}\right) $ with $\left\Vert
P_{i}-Q_{i}\right\Vert <\frac{\epsilon }{n+1}$ $\left( i=1,\cdots ,n\right) $%
, there is a unitary $U\in \mathcal{B}\left( \mathcal{H}\right) $ with $%
\left\Vert U-I\right\Vert <\epsilon $ such that $UP_{i}U^{\ast }=Q_{i}$ for $%
1\leq i\leq n.$
\end{lemma}

\begin{proof}
Define $X=\tsum\limits_{i=1}^{n}Q_{i}P_{i}.$ Let $\delta =\frac{\epsilon }{%
n+1}.$ It is clear that $\sum_{i=1}^{n}P_{i}=\sum_{i=1}Q_{i}=I.$ Since $%
\left\Vert P_{i}-Q_{i}\right\Vert <\delta $ and $P_{i}-Q_{i}$ is
self-adjoint for each $i$, we have that $Q_{i}-P_{i}+\delta \geq 0.$ It
follows that $Q_{i}\geq P_{i}-\delta $ and 
\begin{eqnarray*}
X^{\ast }X &=&\tsum\limits_{i=1}^{n}P_{i}Q_{i}P_{i}\geq
\tsum\limits_{i=1}^{n}P_{i}\left( P_{i}-\delta \right) P_{i} \\
&=&\sum_{i=1}^{n}P_{i}-\sum_{i=1}^{n}\delta P_{i}=\left( 1-\delta \right)
I>0.
\end{eqnarray*}%
Therefore $X$ is invertible and $\left\Vert X^{\ast }X\right\Vert \geq
1-\delta .$ Assume that $X=U\left\vert X\right\vert $ is the polar
decomposition of $X$ where $\left\vert X\right\vert =\left( X^{\ast
}X\right) ^{\frac{1}{2}}$ and $U$ is a partial isometry. Since $X$ is
invertible, we may assume that $U$ is a unitary without loss of generality.
So it is not hard to see that 
\begin{equation*}
\left\Vert \left\vert X\right\vert ^{-1}-I\right\Vert \leq \left( \frac{1}{%
1-\delta }\right) ^{1/2}-1.
\end{equation*}%
Meanwhile, we have $\left\Vert X^{\ast }X\right\Vert \leq 1$ from the
construction of $X$ and the fact that $\left\{ P_{1},\cdots ,P_{n}\right\} $
and $\left\{ Q_{1},\cdots ,Q_{n}\right\} $ are two families of $n$ pairwise
orthogonal projections respectively. Therefore we have that 
\begin{align*}
\left\Vert U-I\right\Vert & \leq \left\Vert U-X\right\Vert +\left\Vert
X-I\right\Vert \\
& \leq \left\Vert X\right\Vert \left\Vert \left\vert X\right\vert
^{-1}-I\right\Vert +\left\Vert \tsum\limits_{i=1}^{n}\left(
Q_{i}-P_{i}\right) P_{i}.\right\Vert \\
& \leq \left( \left( \frac{1}{1-\delta }\right) ^{1/2}-1\right) +n\delta
<\left( n+1\right) \delta =\epsilon .
\end{align*}%
Since $X=\tsum\limits_{i=1}^{n}Q_{i}P_{i}$, it is easy to see $Q_{i}X=XP_{i}$
for $1\leq i\leq n$, then $P_{i}\left\vert X\right\vert =\left\vert
X\right\vert P_{i}$ as well$.$ So 
\begin{equation*}
UP_{i}=X\left\vert X\right\vert ^{-1}P_{i}=XP_{i}\left\vert X\right\vert
^{-1}=Q_{i}X\left\vert X\right\vert ^{-1}=Q_{i}U.
\end{equation*}%
Therefore $UP_{i}U^{\ast }=Q_{i}$ for $1\leq i\leq n$ as desired.
\end{proof}

The following lemma is a useful result concerning the representations of
separable C*-algebras. First, we need to recall that the rank of an operator 
$T\in \mathcal{B}\left( \mathcal{H}\right) ,$ denoted by rank$\left(
T\right) ,$ is the dimension of the closure of the range of $T.$

\begin{lemma}
\label{approx equiv}(Theorem II.5.8., \cite{DK}) Let $\mathcal{A}$ be a
separable unital C*-algebra and $\pi _{i}:\mathcal{A\rightarrow B}\left( 
\mathcal{H}_{i}\right) $ be unital *-representations for $i=1,2.$ Then there
exists a sequence of unitaries $U_{m}:H_{1}\rightarrow H_{2}$ such that $%
\left\Vert \pi _{2}\left( a\right) -U_{m}\pi _{1}\left( a\right) U_{m}^{\ast
}\right\Vert \rightarrow 0$ $\left( m\rightarrow \mathcal{1}\right) $ for
all $a\in \mathcal{A}$ if and only if $rank$ $\left( \pi _{1}\left( a\right)
\right) =rank$ $\left( \pi _{2}\left( a\right) \right) $ for all $a\in 
\mathcal{A}.$
\end{lemma}

\begin{definition}
\label{approx closed}Suppose $\mathcal{H}$ is a separable Hilbert space and $%
F\subseteq \mathcal{H}.$ For given $\epsilon >0,$ we say that%
\begin{equation*}
\left\{ x_{1},\cdots ,x_{n}\right\} \subseteq _{\epsilon }F
\end{equation*}
for $\left\{ x_{1},\cdots ,x_{n}\right\} \subseteq \mathcal{H}$ if there are 
$y_{1},\cdots ,y_{n}\in F$ such that 
\begin{equation*}
\max_{1\leq i\leq n}\left\Vert x_{i}-y_{i}\right\Vert \leq \epsilon .
\end{equation*}
\end{definition}

The following Lemma is a technical result.

\begin{lemma}
\label{20}Let $\mathcal{A\supseteq D\subseteq B}$ be unital inclusions of
separable C*-algebras and $\mathcal{D}$ be a unital finite-dimensional
abelian C*-algebra. Suppose $\mathcal{\rho }_{\mathcal{A}}:\mathcal{%
A\rightarrow B}\left( \mathcal{H}\right) $ and $\rho _{\mathcal{B}}:\mathcal{%
B\rightarrow B}\left( \mathcal{H}\right) $ are representations of $\mathcal{A%
}$ and $\mathcal{B}$ with $\rho _{\mathcal{A}}|_{\mathcal{D}}$ $=\rho _{%
\mathcal{B}}|_{\mathcal{D}}$ on a separable Hilbert space $\mathcal{H}$
respectively. If there are two chains $F_{1}\subseteq F_{2}\subseteq \cdots $
and $G_{1}\subseteq G_{2}\subseteq \cdots $ of finite-dimensional subspaces
of $\mathcal{H}$ satisfying $\dim F_{k}=\dim G_{k}$ for each $k\in \mathbb{N}
$ such that each $F_{k}$ is $\rho _{\mathcal{A}}\left( \mathcal{A}\right) $
invariant and each $G_{k}$ is $\rho _{\mathcal{B}}\left( \mathcal{B}\right) $
invariant$,$ then there are sequences of representations $\left\{ \rho _{k}^{%
\mathcal{A}}\right\} _{k=1}^{\mathcal{1}}$ and $\left\{ \rho _{k}^{\mathcal{B%
}}\right\} _{k=1}^{\mathcal{1}}$ of $\mathcal{A}$ and $\mathcal{B}$ on a
finite-dimensional Hilbert space $\mathcal{H}_{k}$ for each $k$ such that
the restriction of $\widetilde{\rho }_{k}^{\mathcal{A}}$ on subspace $F_{k}$
equals the restriction of $\rho _{\mathcal{A}}$ on\ $F_{k}$, the restriction
of $\widetilde{\rho }_{k}^{\mathcal{B}}$ on subspace $G_{k}$ equals the
restriction of $\rho _{\mathcal{B}}$ on\ $G_{k},$ i.e., 
\begin{equation*}
\widetilde{\rho }_{k}^{\mathcal{A}}|_{F_{k}}=\rho _{\mathcal{A}}|_{F_{k}},%
\text{ \ }\widetilde{\rho }_{k}^{\mathcal{B}}|_{G_{k}}=\rho _{\mathcal{B}%
}|_{G_{k}}
\end{equation*}%
and the restrictions of $\widetilde{\rho }_{k}^{\mathcal{A}}$ and $%
\widetilde{\rho }_{k}^{\mathcal{B}}$ on $\mathcal{D}$ agree, i.e., $%
\widetilde{\rho }_{k}^{\mathcal{A}}|_{\mathcal{D}}=\widetilde{\rho }_{k}^{%
\mathcal{B}}|_{\mathcal{D}}$.
\end{lemma}

\begin{proof}
Suppose that $\mathcal{D}=C^{\ast }\left( p_{1},\cdots ,p_{t}\right) $ where 
$p_{1},\cdots ,p_{t}$ are orthogonal projections with $\tsum%
\limits_{i=1}^{t}p_{i}=I.$ Let $E_{k}=F_{k}+G_{k}$. Note that $E_{k}$ is $%
\rho _{\mathcal{A}}\left( \mathcal{D}\right) $ $\left( =\rho _{\mathcal{B}%
}\left( \mathcal{D}\right) \right) $ invariant. Let $d_{k}=\dim \left(
E_{k}\right) $, $\widetilde{P}_{i}^{k}=\rho _{\mathcal{A}}\left(
p_{i}\right) |_{E_{k}}$ and $r_{i}=rank\left( \widetilde{P}_{i}^{k}\right) .$
Let $E_{k}^{\prime }$ be any finite dimensional subspace of $\mathcal{H}$
that is orthogonal to $E_{k}$ and has dimension $d_{k}^{\prime }=\dim \left(
E_{k}^{\prime }\right) $ so that $d_{k}+d_{k}^{\prime }=l\cdot \dim
F_{k}=l\cdot \dim G_{k}$ and $\frac{rank\left( \rho _{\mathcal{A}}\left(
p_{i}\right) |_{F_{k}}\right) }{\dim \left( F_{k}\right) }\left(
d_{k}+d_{k}^{\prime }\right) =rank\left( \rho _{\mathcal{A}}\left(
p_{i}\right) |_{F_{k}}\right) \cdot l\geq r_{i}$ for $i\in \left\{ 1,\cdots
,t\right\} $, $l\in \mathbb{N}.$ Then we can find projections $\widetilde{Q}%
_{1}^{k},\cdots ,\widetilde{Q}_{t}^{k}\in \mathcal{B}\left( E_{k}^{\prime
}\right) $ such that $\widetilde{Q}_{1}^{k}+\cdots +\widetilde{Q}%
_{t}^{k}=I\in \mathcal{B}\left( E_{k}^{\prime }\right) $, and $%
r_{i}+r_{i}^{^{\prime }}=rank\left( \rho _{\mathcal{A}}\left( p_{i}\right)
|_{F_{k}}\right) \cdot l$ where $r_{i}^{\prime }=rank\left( \widetilde{Q}%
_{i}^{k}\right) .$ Assume that $\mathcal{H}_{k}=E_{k}+E_{k}^{^{\prime }}.$
Since 
\begin{equation*}
\dim \left( \mathcal{H}_{k}\ominus F_{k}\right) =\left( l-1\right) \cdot
\dim F_{k}
\end{equation*}%
and 
\begin{align*}
rank\left( \left( \widetilde{P}_{i}^{k}+\widetilde{Q}_{i}^{k}\right) |_{%
\mathcal{H}_{k}\ominus F_{k}}\right) & =r_{i}+r_{i}^{^{\prime }}-rank\left(
\rho _{\mathcal{A}}\left( p_{i}\right) |_{F_{k}}\right) \\
& =rank\left( \rho _{\mathcal{A}}\left( p_{i}\right) |_{F_{k}}\right) \left(
l-1\right) .
\end{align*}%
We can construct a representation $\rho _{k}^{^{\prime }\mathcal{A}}:%
\mathcal{A\rightarrow B}\left( \mathcal{H}_{k}\ominus F_{k}\right) $ with $%
\rho _{k}^{^{\prime }\mathcal{A}}\left( p_{i}\right) =\left( \widetilde{P}%
_{i}^{k}+\widetilde{Q}_{i}^{k}\right) |_{\mathcal{H}_{k}\ominus F_{k}}$ such
that $\rho _{k}^{\prime \mathcal{A}}$ is unitarily equivalent to the direct
sum of $l-1$ copies of the restriction of $\rho ^{\mathcal{A}}\ $on $F_{k},$%
i.e., $\rho ^{\mathcal{A}}|_{F_{k}}$. Putting $\widetilde{\rho }_{k}^{%
\mathcal{A}}\left( x\right) =\rho ^{\mathcal{A}}\left( x\right)
|_{F_{k}}+\rho _{k}^{^{\prime }\mathcal{A}}\left( x\right) \in \mathcal{B}%
\left( \mathcal{H}_{k}\right) .$ Then $\widetilde{\rho }_{k}^{\mathcal{A}%
}\left( p_{i}\right) =\widetilde{P}_{i}^{k}+\widetilde{Q}_{i}^{k}.$
Similarly, we can construct a representation $\widetilde{\rho }_{k}^{%
\mathcal{B}}$ by the same way such that $\widetilde{\rho }_{k}^{\mathcal{B}%
}\left( p_{i}\right) =\widetilde{P}_{i}^{k}+\widetilde{Q}_{i}^{k}$. This
implies that there are *-representations $\widetilde{\rho }_{k}^{\mathcal{A}%
} $ and $\widetilde{\rho }_{k}^{\mathcal{B}}$ satisfying 
\begin{equation*}
\widetilde{\rho }_{k}^{\mathcal{A}}|_{F_{k}}=\rho _{\mathcal{A}}|_{F_{k}},%
\text{ \ }\widetilde{\rho }_{k}^{\mathcal{B}}|_{G_{k}}=\rho _{\mathcal{B}%
}|_{G_{k}}
\end{equation*}%
and $\widetilde{\rho }_{k}^{\mathcal{A}}|_{\mathcal{D}}=\widetilde{\rho }%
_{k}^{\mathcal{B}}|_{\mathcal{D}}.$
\end{proof}

\bigskip

We need one more technical result for showing our main results in this
section.

\begin{lemma}
\label{21}Let $\mathcal{A\supseteq D\subseteq B}$ be unital inclusions of
C*-algebras in $\tprod\limits_{n=1}^{\mathcal{1}}\mathcal{M}_{k_{n}}\left( 
\mathbb{C}\right) $ and $\mathcal{D}$ be a unital finite-dimensional abelian
C*-subalgebra of $\tprod\limits_{n=1}^{\mathcal{1}}\mathcal{M}_{k_{n}}\left( 
\mathbb{C}\right) .$ Suppose $\Phi :\mathcal{A}\underset{\mathcal{D}}{%
\mathcal{\ast }}\mathcal{B\rightarrow B}\left( \mathcal{H}\right) $ is a
faithful essential representation on a separable Hilbert space $\mathcal{H}.$
Then there are sequences $\left\{ \rho _{m}^{\mathcal{A}}\right\} _{m=1}^{%
\mathcal{1}}$ and $\left\{ \rho _{m}^{\mathcal{B}}\right\} _{m=1}^{\mathcal{1%
}}$ of representations of $\mathcal{A}$ and $\mathcal{B}$ on $\mathcal{H}$
such that $\rho _{m}^{\mathcal{B}}|_{\mathcal{D}}=\rho _{m}^{\mathcal{A}}|_{%
\mathcal{D}}$ and 
\begin{equation*}
\left\Vert \rho _{m}^{\mathcal{A}}\left( a\right) -\Phi _{\mathcal{A}}\left(
a\right) \right\Vert \rightarrow 0\text{ for }\forall a\in \mathcal{A\ }%
\text{as }m\rightarrow \mathcal{1}
\end{equation*}%
\begin{equation*}
\left\Vert \rho _{m}^{\mathcal{B}}\left( b\right) -\Phi _{\mathcal{B}}\left(
b\right) \right\Vert \rightarrow 0\text{ for }\forall b\in \mathcal{B\ }%
\text{as }m\rightarrow \mathcal{1}.
\end{equation*}%
Moreover, for each $m\in \mathbb{N},$ we can find chains of
finite-dimensional subspaces $F_{1}^{m}\subseteq F_{2}^{m}\subseteq \cdots $
and $G_{1}^{m}\subseteq G_{2}^{m}\subseteq \cdots $ of $\mathcal{H}$ with $%
\dim F_{k}^{m}=\dim G_{k}^{m}$ such that each $F_{k}^{m}$ is $\rho _{m}^{%
\mathcal{A}}\left( \mathcal{A}\right) $ invariant, each $G_{k}^{m}$ is $\rho
_{m}^{\mathcal{B}}\left( \mathcal{B}\right) $ invariant and $%
\tbigcup\limits_{k=1}^{\mathcal{1}}F_{k}^{m}$, $\tbigcup\limits_{k=1}^{%
\mathcal{1}}G_{k}^{m}$ are both dense in $\mathcal{H}$.
\end{lemma}

\begin{proof}
Suppose that $\mathcal{D=}C^{\ast }\left( p_{1},\cdots ,p_{t}\right) $ where 
$p_{1},\cdots ,p_{t}$ are orthogonal projections with $\tsum%
\limits_{i=1}^{t}p_{i}=I.$ There are natural *-homomorphisms $\pi _{n}^{%
\mathcal{A}}:\mathcal{A\rightarrow M}_{k_{n}}\left( \mathbb{C}\right) $ and $%
\pi _{n}^{\mathcal{B}}:\mathcal{B\rightarrow M}_{k_{n}}\left( \mathbb{C}%
\right) $ for each $n\in \mathbb{N}$ such that the direct sums of $\left\{
\pi _{n}^{\mathcal{A}}\right\} $ and $\left\{ \pi _{n}^{\mathcal{B}}\right\} 
$ are faithful respectively. We may assume that each $\pi _{k}^{\mathcal{A}}$
and $\pi _{k}^{\mathcal{B}}$ appear infinitely often in the lists $\left\{
\pi _{1}^{\mathcal{A}},\pi _{2}^{\mathcal{A}},\cdots \right\} $ and $\left\{
\pi _{1}^{\mathcal{B}},\pi _{2}^{\mathcal{B}},\cdots \right\} $ respectively
so that we have an increasing sequence $N_{0}=0<N_{1}<N_{2}<\cdots $ such
that $\pi _{k}^{\mathcal{A}}$ and $\pi _{k}^{\mathcal{B}}$ appear at $N_{k}$
position in $\left\{ \pi _{1}^{\mathcal{A}},\pi _{2}^{\mathcal{A}},\cdots
\right\} $ and $\left\{ \pi _{1}^{\mathcal{B}},\pi _{2}^{\mathcal{B}},\cdots
\right\} $ respectively. It is clear that direct sums of them are faithful
representations respectively. Then there are representations $\pi _{\mathcal{%
A}}:$ $\mathcal{A\rightarrow B}\left( \mathcal{H}\right) $ and $\pi _{%
\mathcal{B}}:\mathcal{B\rightarrow B}\left( \mathcal{H}\right) $ with a
projection $P_{N_{k}}$ for each $k\in \mathbb{N}$ such that $P_{N_{k}}$
reduces $\pi _{\mathcal{A}}$ and $\pi _{\mathcal{B}},$ the restrictions of $%
\pi _{\mathcal{A}}$ and $\pi _{\mathcal{B}}$ to $\left(
P_{N_{k}}-P_{N_{k-1}}\right) \left( \mathcal{H}\right) $ are unitarily
equivalent to $\pi _{k}^{\mathcal{A}}$ and $\pi _{k}^{\mathcal{B}}$
respectively, and $P_{N_{k}}\rightarrow I$ in SOT as $k\rightarrow \mathcal{1%
}.$ Since%
\begin{equation*}
rank\text{ }\pi _{\mathcal{A}}\left( a\right) =rank\text{ }\Phi _{\mathcal{A}%
}\left( a\right) \text{ and }rank\text{ }\pi _{\mathcal{B}}\left( b\right)
=rank\text{ }\Phi _{\mathcal{B}}\left( b\right)
\end{equation*}%
for each $a\in \mathcal{A}$ and $b\in \mathcal{B}$ where $\Phi _{\mathcal{A}%
} $ and $\Phi _{\mathcal{B}}$ are the restriction of $\Phi $ on $\mathcal{A}$
and $\mathcal{B}$ respectively, we can find sequences $\left\{ U_{m}\right\}
_{m=1}^{\mathcal{1}}$ and $\left\{ W_{m}\right\} _{m=1}^{\mathcal{1}}$ of
unitaries in $\mathcal{B}\left( \mathcal{H}\right) $ by Lemma \ref{approx
equiv} such that, for every $a\in \mathcal{A}$ and $b\in \mathcal{B},$ we
have%
\begin{equation*}
\left\Vert U_{m}\pi _{\mathcal{A}}\left( a\right) U_{m}^{\ast }-\Phi _{%
\mathcal{A}}\left( a\right) \right\Vert \rightarrow 0\text{ as }m\rightarrow 
\mathcal{1}
\end{equation*}%
\begin{equation*}
\left\Vert W_{m}\pi _{\mathcal{B}}\left( b\right) W_{m}^{\ast }-\Phi _{%
\mathcal{B}}\left( b\right) \right\Vert \rightarrow 0\text{ as }m\rightarrow 
\mathcal{1}.
\end{equation*}%
By the fact that $\Phi _{\mathcal{A}}\left( p_{i}\right) =$ $\Phi _{\mathcal{%
B}}\left( p_{i}\right) $ for every $i\in \left\{ 1,\cdots ,t\right\} ,$ it
follows that 
\begin{equation*}
\left\Vert U_{m}\pi _{\mathcal{A}}\left( p_{i}\right) U_{m}^{\ast }-W_{m}\pi
_{\mathcal{B}}\left( p_{i}\right) W_{m}^{\ast }\right\Vert \rightarrow 0%
\text{ as }\left( m\rightarrow \mathcal{1}\right) .
\end{equation*}%
From Lemma \ref{finite projections}, there is $M_{0}\in \mathbb{N}$ such
that for every $m\geq M_{0},$ there is a unitary $V_{m}$ and $\epsilon _{m}$
satisfying $\left\Vert V_{m}-I\right\Vert <\epsilon _{m}$, $\epsilon
_{m}\rightarrow 0$ $\left( m\rightarrow \mathcal{1}\right) $ and 
\begin{equation*}
V_{m}^{\ast }W_{m}\pi _{\mathcal{B}}\left( p_{i}\right) W_{m}^{\ast
}V_{m}=U_{m}\pi _{\mathcal{A}}\left( p_{i}\right) U_{m}^{\ast }
\end{equation*}%
for each $i\in \left\{ 1,\cdots ,t\right\} $. Without loss of generality we
can assume that, for each $m\in \mathbb{N},$ there is a $V_{m}$ and $%
\epsilon _{m}$ such that $\left\Vert V_{m}-I\right\Vert <\epsilon _{m}$ and 
\begin{equation*}
V_{m}^{\ast }W_{m}\pi _{\mathcal{B}}\left( p_{i}\right) W_{m}^{\ast
}V_{m}=U_{m}\pi _{\mathcal{A}}\left( p_{i}\right) U_{m}^{\ast }
\end{equation*}%
Meanwhile, we still have 
\begin{equation*}
\left\Vert V_{m}^{\ast }W_{m}\pi _{\mathcal{B}}\left( b\right) W_{m}^{\ast
}V_{m}-\Phi _{\mathcal{B}}\left( b\right) \right\Vert \rightarrow 0\text{ }as%
\text{ }m\rightarrow \mathcal{1}.
\end{equation*}%
Let $\rho _{m}^{\mathcal{A}}\left( a\right) =U_{m}\pi _{\mathcal{A}}\left(
a\right) U_{m}^{\ast }$ and $\rho _{m}^{\mathcal{B}}\left( b\right)
=V_{m}^{\ast }W_{m}\pi _{\mathcal{B}}\left( b\right) W_{m}^{\ast }V_{m}$ for
each $m\in \mathbb{N}.$ It is clear that $\rho _{m}^{\mathcal{B}}|_{\mathcal{%
D}}=\rho _{m}^{\mathcal{A}}|_{\mathcal{D}}$ and 
\begin{equation*}
\left\Vert \rho _{m}^{\mathcal{A}}\left( a\right) -\Phi _{\mathcal{A}}\left(
a\right) \right\Vert \rightarrow 0\text{ as }m\rightarrow \mathcal{1}
\end{equation*}%
\begin{equation*}
\left\Vert \rho _{m}^{\mathcal{B}}\left( b\right) -\Phi _{\mathcal{B}}\left(
b\right) \right\Vert \rightarrow 0\text{ as }m\rightarrow \mathcal{1}.
\end{equation*}%
Putting $F_{k}^{m}=U_{m}P_{N_{k}}U_{m}^{\ast }\left( \mathcal{H}\right) $
and $G_{k}^{m}=V_{m}^{\ast }W_{m}P_{N_{k}}W_{m}^{\ast }V_{m}\left( \mathcal{H%
}\right) .$ Note that $\dim F_{k}^{m}=\dim G_{k}^{m}.$ We also have $%
F_{1}^{m}\subseteq F_{2}^{m}\subseteq \cdots $ and $G_{1}^{m}\subseteq
G_{2}^{m}\subseteq \cdots $ are chains of finite dimensional subspaces of $%
\mathcal{H},$ and each $F_{k}^{m}$ is $\rho _{m}^{\mathcal{A}}\left( 
\mathcal{A}\right) $ invariant, each $G_{k}^{m}$ is $\rho _{m}^{\mathcal{B}%
}\left( \mathcal{B}\right) $ invariant. Since $P_{N_{k}}\rightarrow I$ in
SOT as $k\rightarrow \mathcal{1},$ we have $\tbigcup\limits_{k=1}^{\mathcal{1%
}}F_{k}^{m}$ and $\tbigcup\limits_{k=1}^{\mathcal{1}}G_{k}^{m}$ are both
dense in $\mathcal{H}$. This completes the proof.
\end{proof}

From Lemmas \ref{20} and \ref{21}, we are able to obtain the next result
which is a key for proving our main result in this section.

\begin{proposition}
\label{22}Let $\mathcal{A\supseteq D\subseteq B}$ be unital C*-inclusions of
C*-algebras in $\tprod\limits_{n=1}^{\mathcal{1}}\mathcal{M}_{k_{n}}\left( 
\mathbb{C}\right) $ and $\mathcal{D}$ is a unital finite-dimensional abelian
C*-subalgebra. Then $\mathcal{A}\underset{\mathcal{D}}{\mathcal{\ast }}%
\mathcal{B\ }$is RFD.
\end{proposition}

\begin{proof}
Suppose that $\mathcal{D=}C^{\ast }\left( p_{1},\cdots ,p_{t}\right) $ where 
$p_{1},\cdots ,p_{t}$ are orthogonal projections with $\tsum%
\limits_{i=1}^{t}p_{i}=I.$ Suppose $\Phi :\mathcal{A}\underset{\mathcal{D}}{%
\mathcal{\ast }}\mathcal{B\rightarrow B}\left( \mathcal{H}\right) $ is a
faithful essential representation on a separable Hilbert space $\mathcal{H}.$
Then by Lemma \ref{21}, there are sequences $\left\{ \rho _{m}^{\mathcal{A}%
}\right\} $ and $\left\{ \rho _{m}^{\mathcal{B}}\right\} $ of
representations of $\mathcal{A}$ and $\mathcal{B}$ on $\mathcal{H}$ such
that $\rho _{m}^{\mathcal{B}}|_{\mathcal{D}}=\rho _{m}^{\mathcal{A}}|_{%
\mathcal{D}}$ and 
\begin{equation*}
\left\Vert \rho _{m}^{\mathcal{A}}\left( a\right) -\Phi _{\mathcal{A}}\left(
a\right) \right\Vert \rightarrow 0\text{ as }m\rightarrow \mathcal{1}
\end{equation*}%
\begin{equation*}
\left\Vert \rho _{m}^{\mathcal{B}}\left( b\right) -\Phi _{\mathcal{B}}\left(
b\right) \right\Vert \rightarrow 0\text{ as }m\rightarrow \mathcal{1}.
\end{equation*}%
Moreover, for each $m\in \mathbb{N},$ we can find chains of
finite-dimensional subspaces $F_{1}^{m}\subseteq F_{2}^{m}\subseteq \cdots $
and $G_{1}^{m}\subseteq G_{2}^{m}\subseteq \cdots $ of $\mathcal{H}$ with $%
\dim F_{k}^{m}=\dim G_{k}^{m}$ such that each $F_{k}^{m}$ is $\rho _{m}^{%
\mathcal{A}}\left( \mathcal{A}\right) $ invariant, each $G_{k}^{m}$ is $\rho
_{m}^{\mathcal{B}}\left( \mathcal{B}\right) $ invariant, and $%
\tbigcup\limits_{k=1}^{\mathcal{1}}F_{k}^{m}$, $\tbigcup\limits_{k=1}^{%
\mathcal{1}}G_{k}^{m}$ are both dense in $\mathcal{H}$. Then, for each $m\in 
\mathbb{N},$ there are sequences of representations $\left\{ \rho _{m,k}^{%
\mathcal{A}}\right\} _{k=1}^{\mathcal{1}}$ and $\left\{ \rho _{m,k}^{%
\mathcal{B}}\right\} _{k=1}^{\mathcal{1}}$ of $\mathcal{A}$ and $\mathcal{B}$
on a finite-dimensional Hilbert space $\mathcal{H}_{m,k}$ by Lemma \ref{20},
such that 
\begin{equation*}
\widetilde{\rho }_{m,k}^{\mathcal{A}}|_{F_{k}^{m}}=\rho _{m}^{\mathcal{A}%
}|_{F_{k}^{m}},\text{ \ }\widetilde{\rho }_{m,k}^{\mathcal{B}%
}|_{G_{k}^{m}}=\rho _{m}^{\mathcal{B}}|_{G_{k}^{m}}
\end{equation*}%
and $\widetilde{\rho }_{m,k}^{\mathcal{A}}|_{\mathcal{D}}=\widetilde{\rho }%
_{m,k}^{\mathcal{B}}|_{\mathcal{D}}$ for each $k\in \mathbb{N}.$ We first
take representations $\widetilde{\rho }_{1,1}^{\mathcal{B}}$, $\widetilde{%
\rho }_{1,1}^{\mathcal{A}}$ of $\mathcal{A}$ and $\mathcal{B}$ on $\mathcal{H%
}_{1}^{1}$ respectively. Then $\widetilde{\rho }_{1,1}^{\mathcal{B}}\left(
p_{i}\right) =\widetilde{\rho }_{1,1}^{\mathcal{A}}\left( p_{i}\right) $ and 
\begin{equation*}
\widetilde{\rho }_{1,1}^{\mathcal{A}}|_{F_{1}^{1}}=\rho _{1}^{\mathcal{A}%
}|_{F_{1}^{1}},\ \text{ }\widetilde{\rho }_{1,1}^{\mathcal{B}%
}|_{G_{1}^{1}}=\rho _{1}^{\mathcal{B}}|_{G_{1}^{1}}.
\end{equation*}%
Next we can find $F_{l_{2}}^{2}$ and $G_{l_{2}}^{2}$ such that 
\begin{equation*}
\left\{ \eta _{1}^{1},\cdots ,\eta _{t_{1}}^{1}\right\} \subseteq
_{1}G_{l_{2}}^{2}
\end{equation*}%
\begin{equation*}
\left\{ \xi _{1}^{1},\cdots ,\xi _{t_{1}}^{1}\right\} \subseteq
_{1}F_{l_{2}}^{2}
\end{equation*}%
where $\left\{ \xi _{1}^{1},\cdots ,\xi _{t_{1}}^{1}\right\} $ and $\left\{
\eta _{1}^{1},\cdots ,\eta _{t_{1}}^{1}\right\} $ are linear bases of $%
F_{1}^{1}$ and $G_{1}^{1}$ respectively. Moreover, we have representations $%
\widetilde{\rho }_{2,l_{2}}^{\mathcal{A}},\widetilde{\rho }_{2,l_{2}}^{%
\mathcal{B}}$ of $\mathcal{A}$ and $\mathcal{B}$ on $\mathcal{H}_{l_{2}}^{2}$
such that $\widetilde{\rho }_{2,l_{2}}^{\mathcal{B}}\left( p_{i}\right) =%
\widetilde{\rho }_{2,l_{2}}^{\mathcal{A}}\left( p_{i}\right) $ and 
\begin{equation*}
\widetilde{\rho }_{2,l_{2}}^{\mathcal{A}}|_{F_{l_{2}}^{2}}=\rho _{2}^{%
\mathcal{A}}|_{F_{l_{2}}^{2}},\text{ }\widetilde{\rho }_{2,l_{2}}^{\mathcal{B%
}}|_{G_{l_{2}}^{2}}=\rho _{2}^{\mathcal{B}}|_{G_{l_{2}}^{2}}.
\end{equation*}%
Sequentially, we can find $F_{l_{3}}^{3}$ and $G_{l_{3}}^{3}$ satisfying 
\begin{equation*}
\left\{ \xi _{1}^{1},\cdots ,\xi _{t_{1}}^{1},\xi _{1}^{2},\cdots ,\xi
_{t_{m}}^{2}\right\} \subseteq _{\frac{1}{2}}F_{l_{3}}^{3}
\end{equation*}%
\begin{equation*}
\left\{ \eta _{1}^{1},\cdots ,\eta _{t_{1}}^{1},\eta _{1}^{2},\cdots ,\eta
_{t_{m}}^{2}\right\} \subseteq _{\frac{1}{2}}G_{l_{3}}^{3}
\end{equation*}%
where $\left\{ \xi _{1}^{2},\cdots ,\xi _{t_{2}}^{2}\right\} $ and $\left\{
\eta _{1}^{2},\cdots ,\eta _{t_{2}}^{2}\right\} $ are linear bases of $%
F_{l_{2}}^{2}$ and $G_{l_{2}}^{2}$ respectively. Meanwhile, representations $%
\widetilde{\rho }_{3,l_{3}}^{\mathcal{A}},\widetilde{\rho }_{3,l_{3}}^{%
\mathcal{B}}$ of $\mathcal{A}$ and $\mathcal{B}$ are both on $\mathcal{H}%
_{l_{3}}^{3}$ with $\widetilde{\rho }_{3,l_{3}}^{\mathcal{B}}|_{\mathcal{D}}=%
\widetilde{\rho }_{3,l_{3}}^{\mathcal{A}}|_{\mathcal{D}}$ and 
\begin{equation*}
\widetilde{\rho }_{3,l_{3}}^{\mathcal{A}}|_{F_{l_{3}}^{3}}=\rho _{3}^{%
\mathcal{A}}|_{F_{l_{3}}^{3}}\text{, \ \ }\widetilde{\rho }_{3,l_{3}}^{%
\mathcal{B}}|_{G_{l3}^{3}}=\rho _{3}^{\mathcal{B}}|_{G_{l_{3}}^{3}}.
\end{equation*}%
So from the above construction, we can find a sequence $\left\{ \widetilde{%
\rho }_{m,l_{m}}^{\mathcal{B}}\right\} _{m=1}^{\mathcal{1}}$ of
representations and a sequence $\left\{ \widetilde{\rho }_{m,l_{m}}^{%
\mathcal{A}}\right\} _{m}^{\mathcal{1}}$ of representations satisfying $%
\widetilde{\rho }_{m,l_{m}}^{\mathcal{B}}\left( p_{i}\right) =\widetilde{%
\rho }_{m,l_{m}}^{\mathcal{A}}\left( p_{i}\right) $ for each $m\in \mathbb{N}%
.$ We still have that $\tbigcup\limits_{m=1}^{\mathcal{1}}F_{l_{m}}^{m}$, $%
\tbigcup\limits_{m=1}^{\mathcal{1}}G_{l_{m}}^{m}$ are both dense in $%
\mathcal{H}.$ Let $\widetilde{\rho }_{m,l_{m}}:\mathcal{A}\underset{\mathcal{%
D}}{\mathcal{\ast }}\mathcal{B\rightarrow B}\left( \mathcal{H}%
_{l_{m}}^{m}\right) $ be the *-representation such that $\widetilde{\rho }%
_{m,l_{m}}|_{\mathcal{A}}=\widetilde{\rho }_{m,l_{m}}^{\mathcal{A}}$ and $%
\widetilde{\rho }_{m,l_{m}}|_{\mathcal{B}}=\widetilde{\rho }_{m,l_{m}}^{%
\mathcal{B}}.$ We want to show that, for a given $x\in \mathcal{A}\underset{%
\mathcal{D}}{\mathcal{\ast }}\mathcal{B}$ and any $\epsilon >0,$ there is $%
k\in \mathbb{N}$ such that 
\begin{equation*}
\left\Vert \widetilde{\rho }_{k,l_{k}}\left( x\right) \right\Vert \geq
\left\Vert x\right\Vert -\epsilon .
\end{equation*}%
This will suffice to show that $\mathcal{A}\underset{\mathcal{D}}{\mathcal{%
\ast }}\mathcal{B}$ is RFD. Write $x=w_{1}+\cdots +w_{M}$ as the sum of
finitely many words $w_{i}$ in $\mathcal{A}$ and $\mathcal{B}.$ Assume $\xi
\in \mathcal{H}$ is a unit vector such that $\left\Vert \Phi \left( x\right)
\xi \right\Vert \geq \left\Vert \xi \right\Vert -\frac{\epsilon }{2}.$ We
will show that for every $i\in \left\{ 1,\cdots ,M\right\} ,$ there is $%
k\left( i\right) $ such that if $k\geq k\left( i\right) $, then 
\begin{equation*}
\left\Vert \widetilde{\rho }_{k,l_{k}}\left( w_{i}\right) \xi -\Phi \left(
w_{i}\right) \xi \right\Vert <\epsilon /2M.
\end{equation*}%
Taking $k\geq \max_{1\leq i\leq M}k\left( i\right) ,$ this will imply $%
\left\Vert \widetilde{\rho }_{k,l_{k}}\left( x\right) \xi -\Phi \left(
x\right) \xi \right\Vert <\epsilon /2,$ which will yield what we want. To
show it, write 
\begin{equation*}
w_{i}=a_{l}a_{l-1}\cdots a_{2}a_{1}
\end{equation*}%
for some $l\in \mathbb{N}$ and $a_{1},\cdots ,a_{l}\in \mathcal{A\cup B}.$
Let $\xi _{0}=\xi $, $\xi _{j}=\Phi \left( a_{j}\right) \xi _{j-1}$ $\left(
1\leq j\leq l\right) $ and $N=\max_{1\leq j\leq l}\left\Vert
a_{j}\right\Vert .$ Choose $k$ large enough to ensure that 
\begin{equation*}
\max \left( dist\left( \xi _{j-1},F_{l_{k}}^{k}\right) ,dist\left( \xi
_{j-1},G_{l_{k}}^{k}\right) \right) \leq \epsilon /\left( 8lMN^{l-j}\right)
\end{equation*}%
and%
\begin{equation*}
\left\Vert \Phi \left( a_{j}\right) -\rho _{k}^{\mathcal{A}}\left(
a_{j}\right) \right\Vert <\frac{\epsilon }{8lMN^{l-1}}\ \ \ \ \ \ \ \ \text{%
\ if}\ a_{j}\in \mathcal{A}
\end{equation*}%
or 
\begin{equation*}
\left\Vert \Phi \left( a_{j}\right) -\rho _{k}^{\mathcal{B}}\left(
a_{j}\right) \right\Vert <\frac{\epsilon }{8lMN^{l-1}}\ \ \ \ \ \ \ \ \ 
\text{if}\ a_{j}\in \mathcal{B}
\end{equation*}%
for any $j\in \left\{ 1,\cdots ,l\right\} .$ Let $\eta \in \mathcal{H}.$ If $%
a_{j}\in \mathcal{A}$, let $\eta _{k}=P_{F_{l_{k}}^{k}}\left( \eta \right)
\in F_{l_{k}}^{k},$ then 
\begin{align*}
\left\Vert \Phi \left( a_{j}\right) \eta -\widetilde{\rho }_{k,l_{k}}\left(
a_{j}\right) \eta \right\Vert & \leq \left\Vert \Phi \left( a_{j}\right)
\eta -\widetilde{\rho }_{k,l_{k}}\left( a_{j}\right) \eta _{k}\right\Vert
+\left\Vert \widetilde{\rho }_{k,l_{k}}\left( a_{j}\right) \eta _{k}-%
\widetilde{\rho }_{k,l_{k}}\left( a_{j}\right) \eta \right\Vert \\
& \leq \left\Vert \Phi \left( a_{j}\right) \eta -\Phi \left( a_{j}\right)
\eta _{k}+\Phi \left( a_{j}\right) \eta _{k}-\widetilde{\rho }%
_{k,l_{k}}\left( a_{j}\right) \eta _{k}\right\Vert \\
& +\left\Vert \widetilde{\rho }_{k,l_{k}}\left( a_{j}\right) \eta _{k}-%
\widetilde{\rho }_{k,l_{k}}\left( a_{j}\right) \eta \right\Vert \\
& \leq 2\left\Vert a_{j}\right\Vert dist\left( \eta ,F_{l_{k}}^{k}\right)
+\left\Vert \Phi \left( a_{j}\right) \eta _{k}-\rho _{k}^{\mathcal{A}}\left(
a_{j}\right) \eta _{k}\right\Vert \\
& \leq 2\left\Vert a_{j}\right\Vert dist\left( \eta ,F_{l_{k}}^{k}\right) +%
\frac{\epsilon }{4lMN^{l-1}}\left\Vert \eta _{k}\right\Vert
\end{align*}%
Similarly, if $a_{j}\in \mathcal{B},$ then let $\eta
_{k}=P_{G_{l_{k}}^{k}}\left( \eta \right) \in G_{l_{k}}^{k},$ then 
\begin{equation*}
\left\Vert \Phi \left( a_{j}\right) \eta -\widetilde{\rho }_{k,l_{k}}\left(
a_{j}\right) \eta \right\Vert \leq 2\left\Vert a_{j}\right\Vert dist\left(
\eta ,G_{l_{k}}^{k}\right) +\frac{\epsilon }{4lMN^{l-1}}\left\Vert \eta
_{k}\right\Vert
\end{equation*}%
Therefore 
\begin{align*}
\left\Vert \Phi \left( w_{i}\right) \xi -\widetilde{\rho }_{k,l_{k}}\left(
w_{i}\right) \xi \right\Vert & =\left\Vert \Phi \left( a_{l}a_{l-1}\cdots
a_{2}\right) \Phi \left( a_{1}\right) \xi _{0}-\widetilde{\rho }%
_{k,l_{k}}\left( a_{l}a_{l-1}\cdots a_{2}\right) \widetilde{\rho }_{k}\left(
a_{1}\right) \xi _{0}\right\Vert \\
& \leq ||\Phi \left( a_{l}a_{l-1}\cdots a_{2}\right) \Phi \left(
a_{1}\right) \xi _{0}-\widetilde{\rho }_{k,l_{k}}\left( a_{l}a_{l-1}\cdots
a_{2}\right) \Phi \left( a_{1}\right) \xi _{0} \\
& +\widetilde{\rho }_{k}\left( a_{l}a_{l-1}\cdots a_{2}\right) \Phi \left(
a_{1}\right) \xi _{0}-\widetilde{\rho }_{k,l_{k}}\left( a_{l}a_{l-1}\cdots
a_{2}\right) \widetilde{\rho }_{k}\left( a_{1}\right) \xi _{0}|| \\
& \leq \left\Vert \Phi \left( a_{l}a_{l-1}\cdots a_{2}\right) \xi _{1}-%
\widetilde{\rho }_{k,l_{k}}\left( a_{l}a_{l-1}\cdots a_{2}\right) \xi
_{1}\right\Vert \\
& +\left\Vert \widetilde{\rho }_{k,l_{k}}\left( a_{l}a_{l-1}\cdots
a_{2}\right) \right\Vert \left\Vert \Phi \left( a_{1}\right) \xi _{0}-%
\widetilde{\rho }_{k,l_{k}}\left( a_{1}\right) \xi _{0}\right\Vert \\
& \leq \sum_{j=1}^{l-1}\left\Vert \widetilde{\rho }_{k,l_{k}}\left(
a_{l}\cdots a_{j+1}\right) \right\Vert \left\Vert \Phi \left( a_{j}\right)
\xi _{j-1}-\widetilde{\rho }_{k,l_{k}}\left( a_{j}\right) \xi
_{j-1}\right\Vert \\
& <\sum_{j=1}^{l-1}N^{l-j+1}\cdot 2N\cdot \\
& \left( \max \left( dist\left( \xi _{j-1},F_{l_{k}}^{k}\right) ,dist\left(
\xi _{j-1},G_{l_{k}}^{k}\right) \right) +\frac{\epsilon }{4lMN^{l}}%
N^{j-1}\left\Vert \xi _{0}\right\Vert \right) \\
& =\frac{\epsilon }{2M}
\end{align*}%
It follows that $\mathcal{A}\underset{\mathcal{D}}{\mathcal{\ast }}\mathcal{%
B\ }$is RFD.
\end{proof}

The following lemma can be found in [\ref{[BD]}]. Combining previous lemmas
and the lemma below, we are able to state and prove our main result about
unital RFD algebras.

\begin{lemma}
(Lemma 2.2., [\ref{[BD]}]) \label{4}Let $\mathcal{A}$ and $\mathcal{B}$ be
unital C*-algebras having $\mathcal{D}$ embedded as a unital C*-subalgebra
of each of them. Let 
\begin{equation*}
\mathcal{C=A}\underset{\mathcal{D}}{\mathcal{\ast }}\mathcal{B}
\end{equation*}%
be the full amalgamated free product of $\mathcal{A}$ and $\mathcal{B}$ over 
$\mathcal{D}$. If there is a projection $p\in \mathcal{D}$ and there are
partial isometries $v_{1},\cdots ,v_{n}\in \mathcal{D}$ such that $%
v_{i}^{\ast }v_{i}\leq p$ and $\overset{n}{\underset{i=1}{\sum }}%
v_{i}v_{i}^{\ast }=1-p,$ then 
\begin{equation*}
p\mathcal{C}p\cong \left( p\mathcal{A}p\right) \underset{p\mathcal{D}p}{\ast 
}\left( p\mathcal{B}p\right) .
\end{equation*}
\end{lemma}

\begin{remark}
\label{13.3}Suppose $\mathcal{A}$ is a unital C*-algebra and suppose there
is a projection $p\in \mathcal{A}$ and there are partial isometries $%
v_{1},\cdots ,v_{n}\in \mathcal{A}$ such that $v_{i}^{\ast }v_{i}\leq p$ and 
$\overset{n}{\underset{i=1}{\sum }}v_{i}v_{i}^{\ast }=1-p.$ By emulating the
argument in the proof of Lemma 2.1 in [\ref{[BD]}], we know that $\mathcal{A}
$ is MF if and only if $p\mathcal{A}p$ is MF.
\end{remark}

\begin{theorem}
\label{23}Let $\mathcal{A}$, $\mathcal{B}$ be unital RFD algebras and $%
\mathcal{D}$ be a finite-dimensional C*-algebra. Suppose $\psi _{\mathcal{A}%
}:\mathcal{D\rightarrow A}$ and $\psi _{\mathcal{B}}:$ $\mathcal{%
D\rightarrow B}$ are unital embeddings. Then $\mathcal{A}\underset{\mathcal{D%
}}{\mathcal{\ast }}\mathcal{B}$ is RFD if and only if there are unital
embeddings $q_{1}:\mathcal{A\rightarrow }\prod_{n=1}^{\infty }\mathcal{M}%
_{k_{n}}(\mathbb{C})$ and $q_{2}:\mathcal{B\rightarrow }\prod_{n=1}^{\infty }%
\mathcal{M}_{k_{n}}(\mathbb{C})$ for a sequence $\left\{ k_{n}\right\}
_{n=1}^{\mathcal{1}}$ of integers such that the following diagram commutes
\end{theorem}

\begin{equation*}
\begin{array}{ccc}
\mathcal{D} & \overset{\psi _{\mathcal{A}}}{\hookrightarrow } & \mathcal{A}
\\ 
^{\psi _{\mathcal{B}}}\downarrow &  & {\Huge \downarrow }^{q_{1}} \\ 
\mathcal{B} & \overset{q_{2}}{\hookrightarrow } & \prod_{m=1}^{\infty }%
\mathcal{M}_{k_{m}}(\mathbb{C)}%
\end{array}%
\end{equation*}

\begin{proof}
If $\mathcal{A}\underset{\mathcal{D}}{\mathcal{\ast }}\mathcal{B}$ is RFD,
then there is a unital embedding $\Phi :$ $\mathcal{A}\underset{\mathcal{D}}{%
\mathcal{\ast }}\mathcal{B\rightarrow }\prod_{n=1}^{\infty }\mathcal{M}%
_{k_{n}}(\mathbb{C})$ for a sequence $\left\{ k_{n}\right\} _{n=1}^{\infty }$
of integers. Let $q_{1}$ and $q_{2}$ be the restrictions of $\Phi $ on $%
\mathcal{A}$ and $\mathcal{B}$ respectively. Then the above diagram is
commutative. Conversely, we may assume that $\mathcal{A}$, $\mathcal{B}$ are
unital subalgebras of $\prod_{n=1}^{\infty }\mathcal{M}_{k_{n}}(\mathbb{C})$
for a sequence $\left\{ k_{n}\right\} _{n=1}^{\mathcal{1}}$ of integers and $%
\mathcal{A}$ $\supseteq \mathcal{D\subseteq B}$ are unital inclusions of
C*-algebras. Since $\mathcal{D}$ is a finite-dimensional C*-subalgebra, we
can find a projection $p\in \mathcal{D}$ and partial isometries $%
v_{1},\cdots ,v_{n}\in \mathcal{D}$ such that $v_{i}^{\ast }v_{i}\leq p$ and 
$\overset{n}{\underset{i=1}{\sum }}v_{i}v_{i}^{\ast }=1-p.$ Therefore, for
showing $\mathcal{A}\underset{\mathcal{D}}{\mathcal{\ast }}\mathcal{B}$ is
RFD, it is sufficient to show that $P\mathcal{A}P\ast _{P\mathcal{D}P}P%
\mathcal{B}P$ is RFD by Lemma \ref{4} and Lemma 2.1 in \cite{[BD]}. Since $P%
\mathcal{D}P$ is a finite-dimensional abelian C*-algebra. Then the desired
result follows from Proposition \ref{22}.
\end{proof}

\begin{corollary}
Suppose that $\mathcal{A}$ is RFD and $\mathcal{D}$ is a unital
finite-dimensional C*-subalgebra of $\mathcal{A}.$ Then $\mathcal{A\ast }_{%
\mathcal{D}}\mathcal{A}$ is an RFD C*-algebra.
\end{corollary}

\begin{proof}
It is clear by Theorem \ref{23}.
\end{proof}

\begin{example}
Let $\mathcal{M}_{k}\left( \mathbb{C}\right) \supseteq \mathcal{D\subseteq M}%
_{l}\left( \mathbb{C}\right) $ be unital inclusions of unital C*-algebras.
If tr$_{k}|_{\mathcal{D}}=$ tr$_{l}|_{\mathcal{D}}$ where tr$_{k}$ and tr$%
_{l}$ are tracial states on $\mathcal{M}_{k}\left( \mathbb{C}\right) $ and $%
\mathcal{M}_{l}\left( \mathbb{C}\right) $ respectively, then there exists an
integer $n$ and there are two unital embeddings $q_{1}:\mathcal{M}_{k}\left( 
\mathbb{C}\right) \rightarrow \mathcal{M}_{n}\left( \mathbb{C}\right) $ and $%
q_{2}:\mathcal{M}_{l}\left( \mathbb{C}\right) \rightarrow \mathcal{M}%
_{n}\left( \mathbb{C}\right) $ such that $q_{1}|_{\mathcal{D}}=q_{2}|_{%
\mathcal{D}}.$ It implies that there is a commutative diagram which is same
as the one in Theorem \ref{23}. Therefore $\mathcal{M}_{k}\left( \mathbb{C}%
\right) \ast _{\mathcal{D}}\mathcal{M}_{l}\left( \mathbb{C}\right) $ is RFD$%
. $ In fact, this result has been proved in \cite{[BD]}.
\end{example}

\begin{remark}
From the previous example and the fact that every MF algebra has a tracial
state, it is not hard to see that $\mathcal{M}_{k}\left( \mathbb{C}\right)
\ast _{\mathcal{D}}\mathcal{M}_{l}\left( \mathbb{C}\right) $ is RFD if and
only if $\mathcal{M}_{k}\left( \mathbb{C}\right) \ast _{\mathcal{D}}\mathcal{%
M}_{l}\left( \mathbb{C}\right) $ is an MF algebra.
\end{remark}

\section{Full Amalgamated Free Product of Unital MF-Algebras}

\subsection{$\mathcal{D}$ Is A Finite-dimensional C*-algebra}

In this subsection, we consider unital full free products of unital MF
algebras with amalgamation over finite-dimensional C*-algebras. To state and
prove our main result, we need following lemmas.

\begin{lemma}
\label{13.2}Suppose $\mathcal{A=}$C*$\left( x_{1},x_{2},\cdots \right) $ and 
$\mathcal{B=}$C*$\left( y_{1},y_{2},\cdots \right) $ are unital C*-algebras.
Then there is a unital *-homomorphism from $\mathcal{A}$ to $\mathcal{B}$
sending each $x_{k}$ to $y_{k},$ if and only if, for each $\ast $-polynomial 
$P\in \mathbb{C}_{\mathbb{Q}}\left\langle \mathbf{X}_{1},\mathbf{X}%
_{2},\cdots \right\rangle ,$ we have%
\begin{equation*}
\left\Vert P\left( x_{1},x_{2},\cdots \right) \right\Vert \geq \left\Vert
P\left( y_{1},y_{2},\cdots \right) \right\Vert .
\end{equation*}
\end{lemma}

The following lemma is a generalized version of Lemma \ref{approx equiv}.

\begin{lemma}
\label{12}(Theorem 5.1., \cite{HS1}) Suppose $\mathcal{A}$ is a separable
unital C*-algebra, $\mathcal{H}_{1}$, $\mathcal{H}_{2}$ are separable
infinite-dimensional Hilbert spaces and $\pi _{i}:\mathcal{A\rightarrow B(H}%
_{i})$ are unital *-representations for $i=1$, $2.$ If, for each $x\in 
\mathcal{A}$,%
\begin{equation*}
rank(\pi _{1}(x))\leq rank\left( \pi _{2}\left( x\right) \right) ,
\end{equation*}%
then there is a sequence $\left\{ U_{n}\right\} $ of unitary operators, $%
U_{n}:\mathcal{H}_{1}\rightarrow \mathcal{H}_{2},$ such that, for each $x\in 
\mathcal{A}$, 
\begin{equation*}
U_{n}^{\ast }\pi _{2}\left( x\right) U_{n}\rightarrow \pi _{1}\left(
x\right) \text{ }\ast \text{-SOT}
\end{equation*}%
as $n\rightarrow \infty .$
\end{lemma}

The following lemma can be found in \cite{[GH]}, which concerns some
elementary and useful facts about elements in ultraproducts of C*-algebras
and their representatives.

\begin{lemma}
\label{18}(Proposition 2.1, \cite{[GH]}) Let $\mathcal{A}_{i},$ $i\in 
\mathbb{Z},$ be unital C*-algebras and $\alpha $ an ultrafilter on $\mathbb{Z%
}.$ Then

\begin{enumerate}
\item If $P$ is a projection in $\tprod \limits^{\alpha}\mathcal{A}_{l},$
then there are projections $P_{l}$ in $\mathcal{A}_{l}$ such that $P=\left[
\left( P_{l}\right) \right] ;$

\item If $P=\left[ \left( P_{l}\right) \right] ,$ $Q=\left[ \left(
Q_{l}\right) \right] $ are in $\tprod\limits^{\alpha }\mathcal{A}_{l}$ and
all $P_{l},$ $Q_{l}$ are projections and if $V\in \tprod\limits^{\alpha }%
\mathcal{A}_{l}$ is a partial isometry with $V^{\ast }V=P$ and $VV^{\ast
}=Q, $ then there are $V_{l}$ in $\mathcal{A}_{l}$ such that, eventually
along $\alpha ,V_{l}^{\ast }V_{l}=P_{l}$ and $V_{l}V_{l}^{\ast }=Q_{l};$

\item If $P=\left[ \left( P_{l}\right) \right] \in \tprod\limits^{\alpha }%
\mathcal{A}_{l}$ and each $P_{l}$ is a projection, and if $Q$ is a
projection in $\tprod\limits^{\alpha }\mathcal{A}_{l}$ such that $Q\leq P,$
then there are projection $Q_{l}\in \mathcal{A}_{l}$ with $Q_{l}\leq P_{l},$
such that $Q=\left[ \left( Q_{l}\right) \right] .$
\end{enumerate}
\end{lemma}

To show our main result, we need the following technical results.

\begin{lemma}
\label{18.8}Let $\mathcal{A\supset D\subset B}$ be unital inclusions of
MF-algebras. Suppose that $\mathcal{D}$ is a finite-dimensional abelian
C*-algebra generated by a family $\left\{ z_{1},z_{2},\cdots z_{l}\right\} \ 
$of self-adjoint elements$,$ and $\mathcal{A}$ is generated by a family $%
\left\{ x_{1},x_{2}\cdots ,z_{1},z_{2},\cdots z_{l}\right\} \ $of
self-adjoint elements$,\ \mathcal{B}$ is generated by a family $\left\{
y_{1},y_{2},\cdots ,z_{1},z_{2},\cdots z_{l}\right\} \ $of self-adjoint
elements. Let $\varphi :\mathcal{A}\underset{\mathcal{D}}{\mathcal{\ast }}%
\mathcal{B\rightarrow B}\left( \mathcal{H}\right) $ be a faithful
representation of the full amalgamated free product $\mathcal{A}\underset{%
\mathcal{D}}{\mathcal{\ast }}\mathcal{B}$ on a separable Hilbert space $%
\mathcal{H}.$ Assume that there is a sequence $\left\{ k_{n}\right\}
_{n=1}^{\infty }$ of integers with unital embeddings 
\begin{equation*}
q_{1}:\mathcal{A\rightarrow }\prod_{n=1}^{\infty }\mathcal{M}_{k_{n}}(%
\mathbb{C})/\overset{\infty }{\underset{n=1}{\sum }}\mathcal{M}_{k_{n}}(%
\mathbb{C}),
\end{equation*}%
and 
\begin{equation*}
q_{2}:\mathcal{B\rightarrow }\prod_{n=1}^{\infty }\mathcal{M}_{k_{n}}(%
\mathbb{C})/\overset{\infty }{\underset{n=1}{\sum }}\mathcal{M}_{k_{n}}(%
\mathbb{C})
\end{equation*}%
such that $q_{1}\left( z_{i}\right) =q_{2}\left( z_{i}\right) $ for each $%
1\leq i\leq l.$ Also assume that, for a large enough $r\in \mathbb{N},$ 
\begin{equation*}
\left\{ P_{1},\cdots ,P_{2r}\right\} \subset \mathbb{C}_{\mathbb{Q}%
}\left\langle \mathbf{X}_{1},\cdots ,\mathbf{X}_{r},\mathbf{Z}_{1},\cdots ,%
\mathbf{Z}_{l}\right\rangle
\end{equation*}%
with $\left\{ \mathbf{X}_{1},\cdots ,\mathbf{X}_{r},\mathbf{Z}_{1},\cdots ,%
\mathbf{Z}_{l}\right\} \subset \left\{ P_{1},\cdots ,P_{2r}\right\} ,$and 
\begin{equation*}
\left\{ Q_{1},\cdots ,Q_{2r}\right\} \subset \mathbb{C}_{\mathbb{Q}%
}\left\langle \mathbf{Y}_{1},\cdots ,\mathbf{Y}_{r},\mathbf{Z}_{1},\cdots ,%
\mathbf{Z}_{l}\right\rangle
\end{equation*}%
with $\left\{ \mathbf{Y}_{1},\cdots ,\mathbf{Y}_{r},\mathbf{Z}_{1},\cdots ,%
\mathbf{Z}_{l}\right\} \subset \left\{ Q_{1},\cdots ,Q_{2r}\right\} .$ Then
there are sequences $\left\{ E_{m,r}^{i}\right\} _{m=1}^{\mathcal{1}}$, $%
\left\{ F_{m,r}^{i}\right\} _{m=1}^{\mathcal{1}}$ of operators in $\mathcal{B%
}\left( \mathcal{H}\right) $ for each $i\in \mathbb{N},$ and a sequence $%
\left\{ G_{m,r}^{j}\right\} _{m=1}^{\mathcal{1}}$of operators in $\mathcal{B}%
\left( \mathcal{H}\right) $ for each $j\in \left\{ 1,\cdots ,l\right\} $
such that 
\begin{equation}
\left\vert \left\Vert P_{s}\left( E_{m,r}^{1},\cdots
,E_{m,r}^{r},G_{m,r}^{1},\cdots ,G_{m,r}^{l}\right) \right\Vert -\left\Vert
P_{s}\left( x_{1},\cdots ,x_{r},z_{1},\cdots ,z_{l}\right) \right\Vert
\right\vert <\frac{1}{2r},  \tag{5.1}  \label{eq1}
\end{equation}%
\begin{equation}
\left\vert \left\Vert Q_{s}\left( F_{m,r}^{1},\cdots
F_{m,r}^{r},G_{m,r}^{1},\cdots ,G_{m,r}^{l}\right) \right\Vert -\left\Vert
Q_{s}\left( y_{1},\cdots ,y_{r},z_{1},\cdots ,z_{l}\right) \right\Vert
\right\vert <\frac{1}{2r}  \tag{5.2}  \label{eq2}
\end{equation}%
for each $1\leq s\leq 2r,$ $m\in \mathbb{N}.$ We also have that $C^{\ast
}\left( E_{m,r}^{1},\cdots ,F_{m,r}^{1},\cdots ,G_{m,r}^{1},\cdots
,G_{m,r}^{l}\right) $ is an RFD C*-algebra for each $m\in \mathbb{N}$, and 
\begin{align}
E_{m,r}^{i}& \rightarrow \varphi (x_{i})\text{ as }m\rightarrow \infty \text{
in *-SOT for }i\in \mathbb{N};  \tag{5.3}  \label{eq3} \\
F_{m,r}^{i}& \rightarrow \varphi (y_{i})\text{ as }m\rightarrow \infty \text{
in *-SOT for }i\in \mathbb{N};  \tag{5.4}  \label{eq4} \\
G_{m,r}^{j}& \rightarrow \varphi (z_{j})\text{ as }m\rightarrow \infty \text{
in *-SOT for }j\in \left\{ 1,\cdots ,l\right\} .  \tag{5.5}  \label{eq5}
\end{align}
\end{lemma}

\begin{proof}
Without loss of generality, we suppose that $z_{1},\cdots z_{l}$ are
orthogonal projections with $\tsum\limits_{i=1}^{l}z_{i}=I$ and 
\begin{equation*}
\left\Vert x_{i}\right\Vert =\left\Vert y_{i}\right\Vert =1\ \ \ \ \text{%
for\ each}\ i\in \mathbb{N}.
\end{equation*}%
From Lemma \ref{100}, we may assume that for each $i\in \mathbb{N}$, $j\in
\left\{ 1,\cdots ,l\right\} $, there are families 
\begin{equation*}
\left\{ A_{1}^{m},A_{2}^{m},\cdots \right\} ,\left\{
D_{1}^{m},D_{2}^{m},\cdots ,D_{l}^{m}\right\} \text{ and }\left\{
B_{1}^{m},B_{2}^{m},\cdots \right\} \subset \mathcal{M}_{k_{m}}^{s.a}\left( 
\mathbb{C}\right)
\end{equation*}%
for each $k_{m}\in \left\{ k_{n}\right\} _{n=1}^{\mathcal{1}}$ satisfying 
\begin{equation}
\lim_{m\rightarrow \infty }\Vert P(A_{1}^{m},A_{2}^{m}\ldots
,D_{1}^{m},\cdots ,D_{l}^{m})\Vert =\Vert P(x_{1},x_{2}\ldots ,z_{1},\cdots
,z_{l})\Vert  \tag{5.6}  \label{eq6}
\end{equation}%
for any $P\in \mathbb{C}_{\mathbb{Q}}\left\langle \mathbf{X}_{1},\mathbf{X}%
_{2}\cdots ,\mathbf{Z}_{1},\cdots \mathbf{Z}_{l}\right\rangle $, and 
\begin{equation}
\lim_{m\rightarrow \infty }\Vert Q(B_{1}^{m},B_{2}^{m}\ldots
,D_{1}^{m},\cdots ,D_{l}^{m})\Vert =\Vert Q(y_{1},y_{2}\ldots ,z_{1},\cdots
,z_{l})\Vert  \tag{5.7}  \label{eq7}
\end{equation}%
for any $Q\in \mathbb{C}_{\mathbb{Q}}\left\langle \mathbf{Y}_{1},\mathbf{Y}%
_{2},\cdots ,\mathbf{Z}_{1},\cdots \mathbf{Z}_{l}\right\rangle .$

If $r$ is large enough, we can assume that $D_{1}^{m},D_{2}^{m},\cdots
,D_{l}^{m}$ are orthogonal projections with $\sum_{j=1}^{l}D_{j}^{m}=I\in 
\mathcal{M}_{k_{m}}\left( \mathbb{C}\right) $ for $m\geq r$ by Lemma \ref{18}%
. If $\left\{ P_{1},\cdots ,P_{2r}\right\} \subset \mathbb{C}_{\mathbb{Q}%
}\left\langle \mathbf{X}_{1},\cdots ,\mathbf{X}_{r},\mathbf{Z}_{1},\cdots ,%
\mathbf{Z}_{l}\right\rangle $ satisfying 
\begin{equation}
\left\{ \mathbf{X}_{1},\cdots ,\mathbf{X}_{r},\mathbf{Z}_{1},\cdots ,\mathbf{%
Z}_{l}\right\} \subset \left\{ P_{1},\cdots ,P_{2r}\right\}  \tag{5.8}
\label{eq8}
\end{equation}%
and $\left\{ Q_{1},\cdots ,Q_{2r}\right\} \subset \mathbb{C}_{\mathbb{Q}%
}\left\langle \mathbf{Y}_{1},\cdots ,\mathbf{Y}_{r},\mathbf{Z}_{1},\cdots ,%
\mathbf{Z}_{l}\right\rangle $ satisfying 
\begin{equation}
\left\{ \mathbf{Y}_{1},\cdots ,\mathbf{Y}_{r},\mathbf{Z}_{1},\cdots ,\mathbf{%
Z}_{l}\right\} \subset \left\{ Q_{1},\cdots ,Q_{2r}\right\} ,  \tag{5.9}
\label{eq9}
\end{equation}%
then there is an integer $N_{r}$ such that, for $A_{i}\left( N_{r}\right) =%
\underset{r\geqslant N_{r}}{\tprod }A_{i}^{r},B_{i}\left( N_{r}\right) =%
\underset{r\geqslant N_{r}}{\tprod }B_{i}^{r}$ for each $i\in \mathbb{N}$
and $D_{j}\left( N_{r}\right) =\underset{r\geqslant N_{r}}{\tprod }D_{i}^{r}$
for $j\in \left\{ 1,\cdots ,l\right\} ,$ we have 
\begin{equation}
\left\vert \left\Vert P_{s}\left( A_{1}\left( N_{r}\right) ,\cdots
,A_{r}\left( N_{r}\right) ,D_{1}\left( N_{r}\right) ,\cdots ,D_{l}\left(
N_{r}\right) \right) \right\Vert -\left\Vert P_{s}\left( x_{1},\cdots
,x_{r},z_{1},\cdots ,z_{l}\right) \right\Vert \right\vert <\frac{1}{2r} 
\tag{5.10}  \label{eq10}
\end{equation}%
for each $1\leq s\leq 2r$ by (\ref{eq6})$,$ and%
\begin{equation}
\left\vert \left\Vert Q_{t}\left( B_{1}\left( N_{r}\right) ,\cdots
B_{r}\left( N_{r}\right) ,D_{1}\left( N_{r}\right) ,\cdots ,D_{l}\left(
N_{r}\right) \right) \right\Vert -\left\Vert Q_{t}\left( y_{1},\cdots
,y_{r},z_{1},\cdots ,z_{l}\right) \right\Vert \right\vert <\frac{1}{2r} 
\tag{5.11}  \label{eq11}
\end{equation}%
for each $1\leq t\leq 2r$ by (\ref{eq7}). Combining with (\ref{eq8}) and (%
\ref{eq9})$,$ we have 
\begin{equation}
\left\Vert A_{i}(N_{r})\right\Vert \leq 1+\frac{1}{2r}\text{ for }1\leq
i\leq r  \tag{5.12}  \label{eq12}
\end{equation}%
\begin{equation}
\left\Vert B_{i}(N_{r})\right\Vert \leq 1+\frac{1}{2r}for1\leq i\leq r 
\tag{5.13}  \label{eq13}
\end{equation}%
and 
\begin{equation}
\left\Vert D_{i}(N_{r})\right\Vert \leq 1+\frac{1}{2r}\text{ for }1\leq
i\leq l.  \tag{5.14}  \label{eq14}
\end{equation}%
Let 
\begin{equation*}
\mathcal{A}_{N_{r}}=C\text{*}\left( A_{1}(N_{r}),A_{2}(N_{r})\cdots
,D_{1}(N_{r}),\cdots ,D_{l}(N_{r})\right)
\end{equation*}%
and%
\begin{equation*}
\mathcal{B}_{N_{r}}=C\text{*}\left( B_{1}(N_{r}),B_{2}(N_{r})\cdots
,D_{1}(N_{r}),\cdots ,D_{l}(N_{r})\right)
\end{equation*}%
be C*-algebras in $\underset{r\geqslant N_{r}}{\tprod }\mathcal{M}%
_{k_{r}}\left( \mathbb{C}\right) $ and 
\begin{equation*}
\mathcal{D}_{N_{r}}=C\text{*}\left( D_{1}(N_{r}),\cdots ,D_{l}(N_{r})\right)
\end{equation*}%
be unital finite-dimensional C*-subalgebras of $\mathcal{A}_{N_{r}}$ and $%
\mathcal{B}_{N_{r}}$.

Since 
\begin{eqnarray*}
&&\left\Vert P(A_{1}(N_{r}),A_{2}(N_{r})\cdots ,D_{1}(N_{r}),\cdots
,D_{l}(N_{r}))\right\Vert _{\underset{r\geqslant N}{\tprod }\mathcal{M}%
_{k_{r}}\left( \mathbb{C}\right) } \\
&=&\underset{k\geq N_{r}}{\sup }\left\Vert P(A_{1}^{k},A_{2}^{k}\cdots
,D_{1}^{k},\cdots ,D_{l}^{k})\right\Vert _{\mathcal{M}_{k_{r}}\left( \mathbb{%
C}\right) }\geq \left\Vert P(x_{1},x_{2}\cdots ,z_{1},\cdots
,z_{l})\right\Vert _{\mathcal{A}\underset{\mathcal{D}}{\mathcal{\ast }}%
\mathcal{B}}
\end{eqnarray*}%
for any $P\in \mathbb{C}_{\mathbb{Q}}\left\langle \mathbf{X}_{1},\mathbf{X}%
_{2}\cdots ,\mathbf{Z}_{1},\cdots ,\mathbf{Z}_{l}\right\rangle ,$ and 
\begin{align*}
& \left\Vert Q(B_{1}(N_{r}),B_{2}(N_{r})\cdots ,D_{1}(N_{r}),\cdots
,D_{l}(N_{r}))\right\Vert _{\underset{r\geqslant N}{\tprod }\mathcal{M}%
_{k_{r}}\left( \mathbb{C}\right) } \\
& =\underset{k\geq N_{r}}{\sup }\left\Vert Q(B_{1}^{k},B_{2}^{k}\cdots
,D_{1}^{k},\cdots ,D_{l}^{k})\right\Vert _{\mathcal{M}_{k_{r}}\left( \mathbb{%
C}\right) }\geq \left\Vert Q(y_{1},y_{2}\cdots ,z_{1},\cdots
,z_{l})\right\Vert _{\mathcal{A}\underset{\mathcal{D}}{\mathcal{\ast }}%
\mathcal{B}}
\end{align*}%
for any $Q\in \mathbb{C}_{\mathbb{Q}}\left\langle \mathbf{Y}_{1},\mathbf{Y}%
_{2}\cdots ,\mathbf{Z}_{1},\cdots ,\mathbf{Z}_{l}\right\rangle $, there are
*-homomorphisms 
\begin{align*}
\rho _{N_{r}}^{\mathcal{A}}& :\mathcal{A}_{N_{r}}\rightarrow \varphi \left( 
\mathcal{A}\underset{\mathcal{D}}{\mathcal{\ast }}\mathcal{B}\right) ; \\
\rho _{N_{r}}^{\mathcal{B}}& :\mathcal{B}_{Nr}\rightarrow \varphi \left( 
\mathcal{A}\underset{\mathcal{D}}{\mathcal{\ast }}\mathcal{B}\right) .
\end{align*}%
such that $\rho _{N_{r}}^{\mathcal{A}}\left( A_{i}\left( N_{r}\right)
\right) =\varphi \left( x_{i}\right) ,\ \rho _{N_{r}}^{\mathcal{B}}\left(
B_{i}\left( N_{r}\right) \right) =\varphi \left( y_{i}\right) $ and $\rho
_{N_{r}}^{\mathcal{A}}\left( D_{j}\left( N_{r}\right) \right) =\rho
_{N_{r}}^{\mathcal{B}}\left( D_{j}\left( N_{r}\right) \right) =\varphi
\left( z_{j}\right) $ for $i\in \mathbb{N}$ and $j\in \left\{ 1,\cdots
,l\right\} $ by Lemma \ref{13.2}. It follows that there is a *-homomorphism%
\begin{equation*}
\rho _{N_{r}}:\mathcal{A}_{N_{r}}\underset{_{\mathcal{D}_{N_{r}}}}{\ast }%
\mathcal{B}_{N_{r}}\rightarrow \varphi \left( \mathcal{A}\underset{\mathcal{D%
}}{\mathcal{\ast }}\mathcal{B}\right)
\end{equation*}%
satisfying $\rho _{N_{r}}\left( A_{i}\left( N_{r}\right) \right) =\varphi
\left( x_{i}\right) $ and $\rho _{N_{r}}\left( B_{i}\left( N_{r}\right)
\right) =\varphi \left( y_{i}\right) $ as well as $\rho _{N_{r}}\left(
D_{j}\left( N_{r}\right) \right) =\varphi \left( z_{j}\right) $ for each $%
i\in \mathbb{N}\ $and $j\in \left\{ 1,\cdots ,l\right\} .$ We also know that 
$\mathcal{A}_{N_{r}}\underset{_{\mathcal{D}_{N_{r}}}}{\ast }\mathcal{B}%
_{N_{r}}$ is an RFD C*-algebra by Theorem \ref{23}.

Let $\pi _{N_{r}}:\mathcal{A}_{N_{r}}\underset{_{\mathcal{D}_{N_{r}}}}{\ast }%
\mathcal{B}_{N_{r}}\rightarrow \mathcal{B}\left( \mathcal{H}_{N_{r}}\right) $
be a faithful essential representation of $\mathcal{A}_{N_{r}}\underset{_{%
\mathcal{D}_{N_{r}}}}{\ast }\mathcal{B}_{N_{r}}.$ Then $\pi _{N_{r}}\left( 
\mathcal{A}_{N_{r}}\underset{_{\mathcal{D}_{N_{r}}}}{\ast }\mathcal{B}%
_{N_{r}}\right) $ is an RFD C*-algebra and 
\begin{equation}
\pi _{N_{r}}\left( \mathcal{A}_{N_{r}}\underset{_{\mathcal{D}_{N_{r}}}}{\ast 
}\mathcal{B}_{N_{r}}\right) =C^{\ast }\left( \pi _{N_{r}}\left( A_{1}\left(
N_{r}\right) \right) ,\cdots ,\pi _{N_{r}}\left( B_{1}\left( N_{r}\right)
\right) ,\cdots ,D_{1}\left( N_{r}\right) ,\cdots ,D_{l}\left( N_{r}\right)
\right) .  \tag{5.15}  \label{eq15}
\end{equation}%
Since 
\begin{equation*}
rank\left( \pi _{N_{r}}(x)\right) \geq rank\left( \rho _{N_{r}}\left(
x\right) \right)
\end{equation*}%
for every $x\in \mathcal{A}_{N_{r}}\underset{_{\mathcal{D}_{N_{r}}}}{\ast }%
\mathcal{B}_{N_{r}}$, Lemma \ref{12} implies that there is a sequence of
unitary operators $\left\{ U_{m}^{N_{r}}\right\} _{m=1}^{\infty }\subset 
\mathcal{B(H},\mathcal{H}_{N_{r}}\mathcal{)}$ such that%
\begin{equation}
\rho _{N_{r}}\left( x\right) =\ast \text{-}SOT\text{-}\underset{m\rightarrow
\infty }{\lim }U_{m}^{N_{r}\ast }\pi _{N_{r}}\left( x\right) U_{m}^{N_{r}} 
\tag{5.16}  \label{eq16}
\end{equation}%
for each $x\in \mathcal{A}_{N_{r}}\underset{_{\mathcal{D}_{N_{r}}}}{\ast }%
\mathcal{B}_{N_{r}}$. So, for $i\in \mathbb{N}$ and $j\in \left\{ 1,\cdots
,l\right\} ,$ if we put 
\begin{equation}
E_{m,r}^{i}=U_{m}^{N_{r}\ast }\pi _{N_{r}}\left( A_{i}\left( N_{r}\right)
\right) U_{m}^{N_{r}},  \tag{5.17}  \label{eq17}
\end{equation}%
$\ $%
\begin{equation}
F_{m,r}^{i}=U_{m}^{N_{r}\ast }\pi _{N_{r}}\left( B_{i}\left( N_{r}\right)
\right) U_{m}^{N_{r}}  \tag{5.18}  \label{eq18}
\end{equation}%
and 
\begin{equation}
G_{m,r}^{j}=U_{m}^{N_{r}\ast }\pi _{N_{r}}\left( D_{j}\left( N_{r}\right)
\right) U_{m}^{N_{r}},  \tag{5.19}  \label{eq19}
\end{equation}%
then, for every $m\in \mathbb{N},i\in \mathbb{N},j\in \left\{ 1,\cdots
,l\right\} $, we have 
\begin{equation}
\left\Vert P\left( E_{m,r}^{1},\cdots ,E_{m,r}^{r},G_{m,r}^{1},\cdots
,G_{m,r}^{l}\right) \right\Vert =\left\Vert P\left( A_{1}\left( N_{r}\right)
,\cdots ,A_{r}\left( N_{r}\right) ,D_{1}\left( N_{r}\right) ,\cdots
,D_{l}\left( N_{r}\right) \right) \right\Vert ,  \tag{5.20}  \label{eq20}
\end{equation}%
for every $P\in \mathbb{C}\left\langle X_{1},X_{2},\cdots ,Z_{1},\cdots
,Z_{l}\right\rangle ,$ and%
\begin{equation}
\left\Vert Q\left( F_{m,r}^{1},\cdots F_{m,r}^{r},G_{m,r}^{1},\cdots
,G_{m,r}^{l}\right) \right\Vert =\left\Vert Q\left( B_{1}\left( N_{r}\right)
,\cdots B_{r}\left( N_{r}\right) ,D_{1}\left( N_{r}\right) ,\cdots
,D_{l}\left( N_{r}\right) \right) \right\Vert  \tag{5.21}  \label{eq21}
\end{equation}%
for every $Q\in \mathbb{C}\left\langle Y_{1},Y_{2},\cdots ,Z_{1},\cdots
,Z_{l}\right\rangle .$ It follows that, for each $P_{s}\in \left\{
P_{1},\cdots ,P_{2r}\right\} \subseteq \mathbb{C}_{\mathbb{Q}}\left\langle
X_{1},\cdots ,X_{r},Z_{1},\cdots ,Z_{l}\right\rangle $ and each $Q_{t}\in
\left\{ Q_{1},\cdots ,Q_{2r}\right\} \subseteq \mathbb{C}_{\mathbb{Q}%
}\left\langle Y_{1},\cdots ,Y_{r},Z_{1},\cdots ,Z_{l}\right\rangle $ 
\begin{equation*}
\left\vert \left\Vert P_{s}\left( E_{m,r}^{1},\cdots
,E_{m,r}^{r},G_{m,r}^{1},\cdots ,G_{m,r}^{l}\right) \right\Vert -\left\Vert
P_{s}\left( x_{1},\cdots ,x_{r},z_{1},\cdots ,z_{l}\right) \right\Vert
\right\vert <\frac{1}{2r},
\end{equation*}%
\begin{equation*}
\left\vert \left\Vert Q_{t}\left( F_{m,r}^{1},\cdots
F_{m,r}^{r},G_{m,r}^{1},\cdots ,G_{m,r}^{l}\right) \right\Vert -\left\Vert
Q_{t}\left( y_{1},\cdots ,y_{r},z_{1},\cdots ,z_{l}\right) \right\Vert
\right\vert <\frac{1}{2r}
\end{equation*}%
for each $m\in \mathbb{N}$ by (\ref{eq10}), (\ref{eq11}), and (\ref{eq20}), (%
\ref{eq21}). Since $C^{\ast }\left( E_{m,r}^{1},\cdots ,F_{m,r}^{1},\cdots
,G_{m,r}^{1},\cdots ,G_{m,r}^{l}\right) $ is *-isomorphic to the C*-algebra 
\begin{equation*}
C^{\ast }\left( \pi _{N_{r}}\left( A_{1}\left( N_{r}\right) \right) ,\cdots
,\pi _{N_{r}}\left( B_{1}\left( N_{r}\right) \right) ,\cdots ,D_{1}\left(
N_{r}\right) ,\cdots ,D_{l}\left( N_{r}\right) \right) ,
\end{equation*}%
we have 
\begin{equation*}
C^{\ast }\left( E_{m,r}^{1},\cdots ,F_{m,r}^{1},\cdots ,G_{m,r}^{1},\cdots
,G_{m,r}^{l}\right)
\end{equation*}%
is an RFD C*-algebra for each $m\in \mathbb{N}$. We also get 
\begin{align*}
E_{m,r}^{i}& \rightarrow \varphi (x_{i})\text{ as }m\rightarrow \infty \text{
in *-SOT for }i\in \mathbb{N}; \\
F_{m,r}^{i}& \rightarrow \varphi (y_{i})\text{ as }m\rightarrow \infty \text{
in *-SOT for }i\in \mathbb{N}; \\
G_{m,r}^{j}& \rightarrow \varphi (z_{j})\text{ as }m\rightarrow \infty \text{
in *-SOT for }j\in \left\{ 1,\cdots ,l\right\} .
\end{align*}%
by (\ref{eq16}), (\ref{eq17}), (\ref{eq18}) and (\ref{eq19}).
\end{proof}

The next proposition is a key ingredient for proving our main theorem in
this subsection.

\begin{proposition}
\label{15}Let $\mathcal{A\supset D\subset B}$ be unital inclusions of
MF-algebras, where $\mathcal{D}$ is a finite-dimensional abelian C*-algebra.
Let $\varphi :\mathcal{A}\underset{\mathcal{D}}{\mathcal{\ast }}\mathcal{%
B\rightarrow B}\left( \mathcal{H}\right) $ be a faithful representation of
the full amalgamated free product $\mathcal{A}\underset{\mathcal{D}}{%
\mathcal{\ast }}\mathcal{B}$. Suppose that $\mathcal{D}$ is generated by a
family $\left\{ z_{1},z_{2},\cdots z_{l}\right\} \ $of self-adjoint elements$%
,\ \mathcal{A}$ is generated by a family $\left\{ x_{1},x_{2}\cdots
,z_{1},z_{2},\cdots z_{l}\right\} $ of self-adjoint elements and $\mathcal{B}
$ is generated by a family $\left\{ y_{1},y_{2},\cdots ,z_{1},z_{2},\cdots
z_{l}\right\} $ of self-adjoint elements. Suppose that there is a sequence $%
\left\{ k_{n}\right\} _{n=1}^{\infty }$ of integers with unital embeddings $%
q_{1}:\mathcal{A\rightarrow }\prod_{n=1}^{\infty }\mathcal{M}_{k_{n}}(%
\mathbb{C})/\overset{\infty }{\underset{n=1}{\sum }}\mathcal{M}_{k_{n}}(%
\mathbb{C}),$ and $q_{2}:\mathcal{B\rightarrow }\prod_{n=1}^{\infty }%
\mathcal{M}_{k_{n}}(\mathbb{C})/\overset{\infty }{\underset{n=1}{\sum }}%
\mathcal{M}_{k_{n}}(\mathbb{C})$ such that $q_{1}\left( z_{i}\right)
=q_{2}\left( z_{i}\right) $ for each $1\leq i\leq l.$ Then there is a
sequence $\left\{ t_{m}\right\} _{m=1}^{\mathcal{1}}$ of integers such that,
for each $t_{r}\in \left\{ t_{m}\right\} _{m=1}^{\mathcal{1}},$ there exist
sequences 
\begin{equation*}
\left\{ X_{1}^{r},X_{2}^{r},\cdots \right\} ,\left\{
Y_{1}^{r},Y_{2}^{r},\cdots \right\} \text{ and }\left\{ Z_{1}^{r},\cdots
Z_{l}^{r}\right\}
\end{equation*}%
in $\mathcal{M}_{t_{r}}\left( \mathbb{C}\right) $ and a unitary operator $%
W_{r}:\mathcal{H\rightarrow }\left( \mathbb{C}^{t_{r}}\right) ^{\infty }$
satisfying 
\begin{equation*}
W_{r}^{\ast }\left( X_{i}^{r}\right) ^{\left( \infty \right)
}W_{r}\rightarrow \varphi \left( x_{i}\right) \text{ in SOT as }r\rightarrow
\infty \text{ for }i\in \mathbb{N}
\end{equation*}%
\begin{equation*}
W_{r}^{\ast }\left( Y_{i}^{r}\right) ^{\left( \infty \right)
}W_{r}\rightarrow \varphi \left( y_{i}\right) \text{ in SOT as }r\rightarrow
\infty \text{ for }i\in \mathbb{N}
\end{equation*}%
and%
\begin{equation*}
W_{r}^{\ast }\left( Z_{i}^{r}\right) ^{\left( \infty \right)
}W_{r}\rightarrow \varphi \left( z_{i}\right) \text{ in SOT as }r\rightarrow
\infty \text{ for }i\in \left\{ 1,\cdots ,l\right\}
\end{equation*}%
as well as 
\begin{align*}
\left\Vert P(x_{1},x_{2},\cdots ,z_{1},\cdots z_{l})\right\Vert &
=\lim_{r\rightarrow \infty }\left\Vert P\left( X_{1}^{r},X_{2}^{r},\cdots
,Z_{1}^{r},\cdots ,Z_{l}^{r}\right) \right\Vert \\
\left\Vert Q(y_{1},y_{2},\cdots ,z_{1},\cdots z_{l})\right\Vert &
=\lim_{r\rightarrow \infty }\left\Vert Q\left( Y_{1}^{r},Y_{2}^{r},\cdots
,Z_{1}^{r},\cdots ,Z_{l}^{r}\right) \right\Vert
\end{align*}%
for any $P\in \mathbb{C}_{\mathbb{Q}}\left\langle \mathbf{X}_{1},\mathbf{X}%
_{2}\cdots ,\mathbf{Z}_{1},\cdots \mathbf{Z}_{l}\right\rangle $ and $Q\in 
\mathbb{C}_{\mathbb{Q}}\left\langle \mathbf{Y}_{1},\mathbf{Y}_{2}\cdots ,%
\mathbf{Z}_{1},\cdots \mathbf{Z}_{l}\right\rangle $.
\end{proposition}

\begin{proof}
Suppose $z_{1},\cdots z_{l}$ are orthogonal projections with $%
\tsum\limits_{i=1}^{l}z_{i}=I$ and 
\begin{equation*}
\left\Vert x_{i}\right\Vert =\left\Vert y_{i}\right\Vert =1\ \ \ \text{\
for\ each}\ i\in \mathbb{N}.
\end{equation*}%
Assume $\left\{ e_{1},e_{2},\cdots \right\} $ is a family of orthonormal
basis of $\mathcal{H}.$ With notations as in Lemma \ref{18.8}, for a large
enough integer $r$ and a subset $\left\{ e_{1},\cdots ,e_{r}\right\}
\subseteq \left\{ e_{1},e_{2}\cdots \right\} $, there is an integer $M$ such
that, for $1\leq i\leq r,j\in \left\{ 1,\cdots ,l\right\} $ and $1\leq k\leq
r$%
\begin{equation}
\left\Vert E_{M,r}^{i}e_{k}-\varphi (x_{i})e_{k}\right\Vert <\frac{1}{2r} 
\tag{5.22}  \label{eq22}
\end{equation}%
\begin{equation}
\left\Vert F_{M,r}^{i}e_{k}-\varphi (y_{j})e_{k}\right\Vert <\frac{1}{2r} 
\tag{5.23}  \label{eq23}
\end{equation}%
and%
\begin{equation}
\left\Vert G_{M,r}^{j}e_{k}-\varphi (z_{j})e_{k}\right\Vert <\frac{1}{2r}. 
\tag{5.24}  \label{eq24}
\end{equation}

Note that $\left\{ E_{M,r}^{1},E_{M,r}^{2},\cdots ,F_{M,r}^{1},\cdots
,G_{M,r}^{1},\cdots G_{M,r}^{l}\right\} $ is a family of self-adjoint
elements in $\mathcal{B}\left( \mathcal{H}\right) $ from the proof of Lemma %
\ref{18.8}. So, by Lemma \ref{10} and the fact that $C^{\ast }\left(
E_{M,r}^{1},E_{M,r}^{2},\cdots ,F_{M,r}^{1},\cdots ,G_{M,r}^{1},\cdots
G_{M,r}^{l}\right) $ is a QD algebra (actually it is RFD C*-algebra), there
is a projection $\mathcal{P}_{r}\in \mathcal{B(H)}$ such that, for $1\leq
i\leq r,$ $j\in \left\{ 1,\cdots ,l\right\} $, $1\leq k\leq r,$%
\begin{equation}
\left\Vert e_{k}-\mathcal{P}_{r}e_{k}\right\Vert <\frac{1}{6r}  \tag{5.25}
\label{eq25}
\end{equation}%
and%
\begin{eqnarray}
\left\Vert \mathcal{P}_{r}E_{M,r}^{i}\mathcal{P}_{r}e_{k}-E_{M,r}^{i}e_{k}%
\right\Vert &<&\frac{1}{6r}  \TCItag{5.26}  \label{eq26} \\
\left\Vert \mathcal{P}_{r}F_{M,r}^{i}\mathcal{P}_{r}e_{k}-F_{M,r}^{i}e_{k}%
\right\Vert &<&\frac{1}{6r}  \TCItag{5.27}  \label{eq27} \\
\left\Vert \mathcal{P}_{r}G_{M,r}^{j}\mathcal{P}_{r}e_{k}-G_{M,r}^{j}e_{k}%
\right\Vert &<&\frac{1}{6r}  \TCItag{5.28}  \label{eq28}
\end{eqnarray}%
as well as 
\begin{align}
& |\left\Vert P_{s}\left( \mathcal{P}_{r}E_{M,r}^{1}\mathcal{P}_{r},\cdots ,%
\mathcal{P}_{r}E_{M,r}^{r}\mathcal{P}_{r},\mathcal{P}_{r}G_{M,r}^{1}\mathcal{%
P}_{r},\cdots \mathcal{P}_{r}G_{M,r}^{l}\mathcal{P}_{r}\right) \right\Vert
-\left\Vert P_{s}\left( E_{M,r}^{1},\cdots E_{M,r}^{r},G_{M,r}^{1},\cdots
G_{M,r}^{l}\right) \right\Vert |  \notag \\
& <\frac{1}{2r}\text{ \qquad for }1\leq s\leq 2r,\qquad \qquad \qquad \qquad 
\text{\qquad\ \ \ \ \ \ \ \ \ \ \ \ \ \ \ \ \ \ \ \ \ \ \ \ \ \ \ \ }\qquad 
\tag{5.29}  \label{eq29}
\end{align}%
\begin{align}
& |\left\Vert Q_{t}\left( \mathcal{P}_{r}F_{M,r}^{1}\mathcal{P}_{r},\cdots ,%
\mathcal{P}_{r}F_{M,r}^{r}\mathcal{P}_{r},\mathcal{P}_{r}G_{M,r}^{1}\mathcal{%
P}_{r},\cdots \mathcal{P}_{r}G_{M,r}^{l}\mathcal{P}_{r}\right) \right\Vert
-\left\Vert Q_{t}\left( F_{M,r}^{1},\cdots F_{M,r}^{r},G_{M,r}^{1},\cdots
G_{M,r}^{l}\right) \right\Vert |  \notag \\
& <\frac{1}{2r}\text{ \qquad \qquad for }1\leq t\leq 2r.\qquad \qquad \qquad
\qquad \ \ \ \ \   \tag{5.30}  \label{eq30}
\end{align}%
By (\ref{eq20}), (\ref{eq21}) and (\ref{eq29}), (\ref{eq30}), we have that%
\begin{align}
\ & |\left\Vert P_{s}\left( \mathcal{P}_{r}E_{M,r}^{1}\mathcal{P}_{r},\cdots
,\mathcal{P}_{r}E_{M,r}^{r}\mathcal{P}_{r},\mathcal{P}_{r}G_{M,r}^{1}%
\mathcal{P}_{r},\cdots \mathcal{P}_{r}G_{M,r}^{l}\mathcal{P}_{r}\right)
\right\Vert  \notag \\
& -\left\Vert P_{s}\left( A_{1}\left( N_{r}\right) ,\cdots ,A_{r}\left(
N_{r}\right) ,C_{1}\left( N_{r}\right) ,\cdots C_{l}\left( N_{r}\right)
\right) \right\Vert |  \notag \\
& <\frac{1}{2r}\text{ for }1\leq s\leq 2r\ \ \ \ \ \ \ \ \ \ \ \ \ \ \ \ \ \
\   \tag{5.31}  \label{eq31}
\end{align}%
\begin{align}
\ \ & |\left\Vert Q_{t}\left( \mathcal{P}_{r}F_{M,r}^{1}\mathcal{P}%
_{r},\cdots ,\mathcal{P}_{r}F_{M,r}^{r}\mathcal{P}_{r},\mathcal{P}%
_{r}G_{M,r}^{1}\mathcal{P}_{r},\cdots \mathcal{P}_{r}G_{M,r}^{l}\mathcal{P}%
_{r}\right) \right\Vert  \notag \\
-& \left\Vert Q_{t}\left( B_{1}\left( N_{r}\right) ,\cdots B_{r}\left(
N_{r}\right) ,C_{1}\left( N_{r}\right) ,\cdots C_{l}\left( N_{r}\right)
\right) \right\Vert |  \notag \\
& <\frac{1}{2r}\text{ \qquad \qquad for }1\leq t\leq 2r.\qquad \qquad \qquad
\qquad \ \ \ \ \   \tag{5.32}  \label{eq32}
\end{align}%
Let $t_{r}=\dim \mathcal{P}_{r}\mathcal{H}$ and $\widetilde{W}_{r}:\mathcal{P%
}_{r}\mathcal{H\rightarrow }\mathbb{C}^{t_{r}}$ be an unitary$.$ Putting $%
X_{i}^{r}=\widetilde{W}_{r}\mathcal{P}_{r}E_{M,r}^{i}\mathcal{P}_{r}%
\widetilde{W}_{r}^{\ast }$, $Y_{i}^{r}=\widetilde{W}_{r}\mathcal{P}%
_{r}F_{M,r}^{i}\mathcal{P}_{r}\widetilde{W}_{r}^{\ast }$ for $i\in \mathbb{N}
$ and $Z_{j}^{r}=\widetilde{W}_{r}\mathcal{P}_{r}G_{M,r}^{j}\mathcal{P}_{r}%
\widetilde{W}_{r}^{\ast }$ for $j\in \left\{ 1,\cdots ,l\right\} ,$ and
combining (\ref{eq31}), (\ref{eq32}) and (\ref{eq10}), (\ref{eq11})$,$ we
have%
\begin{align}
\left\vert \left\Vert P_{s}\left( X_{1}^{r},\cdots
,X_{r}^{r},Z_{1}^{r},\cdots ,Z_{l}^{r}\right) \right\Vert -\left\Vert
P_{s}(x_{1},\cdots ,x_{r},z_{1},\cdots ,z_{l}\right\Vert \right\vert & <%
\frac{1}{r}\text{ }  \tag{5.33}  \label{eq33} \\
\left\vert \left\Vert Q_{t}\left( Y_{1}^{r},\cdots
,Y_{r}^{r},Z_{1}^{r},\cdots ,Z_{l}^{r}\right) \right\Vert -\left\Vert
Q_{t}(y_{1},\cdots ,y_{r},z_{1},\cdots ,z_{l}\right\Vert \right\vert & <%
\frac{1}{r}\text{ }  \tag{5.34}  \label{eq34}
\end{align}%
for $1\leq s,t\leq 2r.$ Hence we can find a unitary $W_{r}$ $:\mathcal{H}%
\rightarrow \left( \mathbb{C}^{t_{r}}\right) ^{\mathcal{1}}$ such that $W_{r}%
\mathcal{\ }$is unitary equivalent to $\left( \widetilde{W}_{r}\right) ^{%
\mathcal{1}}$ and $W_{r}\mathcal{P}_{r}=\widetilde{W}_{r}$. It follows that,
for $1\leq i\leq r$, $j\in \left\{ 1,\cdots ,l\right\} $ and $1\leq k\leq r,$
we have 
\begin{align}
& \left\Vert W_{r}^{\ast }\left( X_{i}^{r}\right) ^{\infty
}W_{r}e_{k}-E_{m_{r}}^{i}e_{k}\right\Vert \leq \left\Vert W_{r}^{\ast
}\left( X_{i}^{r}\right) ^{\infty }W_{r}\right\Vert \left\Vert e_{k}-%
\mathcal{P}_{r}e_{k}\right\Vert +\left\Vert W_{r}^{\ast }\left(
X_{i}^{r}\right) ^{\infty }W_{r}\mathcal{P}_{r}e_{k}-E_{m_{r}}^{i}e_{k}%
\right\Vert  \notag \\
& \leq \left( 1+\frac{1}{2r}\right) \frac{1}{6r}+\left\Vert W_{r}^{\ast
}\left( \widetilde{W}_{r}\mathcal{P}_{r}E_{M,r}^{i}\mathcal{P}_{r}\widetilde{%
W}_{r}^{\ast }\right) ^{\infty }\widetilde{W}_{r}e_{k}-E_{m_{r}}^{i}e_{k}%
\right\Vert  \notag \\
& =\left( 1+\frac{1}{2r}\right) \frac{1}{6r}+\left\Vert \mathcal{P}%
_{r}E_{M,r}^{i}\mathcal{P}_{r}e_{k}-E_{M,r}^{i}e_{k}\right\Vert  \notag \\
& \leq \left( 1+\frac{1}{2r}\right) \frac{1}{6r}+\frac{1}{6r}<\frac{1}{2r}%
\text{,\qquad }  \tag{5.35}  \label{eq35} \\
\text{and }& \left\Vert W_{r}^{\ast }\left( Y_{i}^{r}\right) ^{\infty
}W_{r}e_{k}-F_{M,r}^{i}e_{k}\right\Vert <\frac{1}{2r}\text{,\qquad } 
\tag{5.36}  \label{eq36} \\
\ & \left\Vert W_{r}^{\ast }\left( Z_{j}^{r}\right) ^{\infty
}W_{r}e_{k}-G_{M,r}^{j}e_{k}\right\Vert <\frac{1}{2r}  \tag{5.37}
\label{eq37}
\end{align}%
by the definition of $X_{i}^{r}$, $Y_{i}^{r}$ and $Z_{j}^{r}$ and (\ref{eq25}%
), (\ref{eq26}), (\ref{eq27}) and (\ref{eq28})$.$ Combining the inequalities
from above with (\ref{eq22}), (\ref{eq23}) and (\ref{eq24}), we have, for $%
1\leq i\leq r$, $j\in \left\{ 1,\cdots ,l\right\} $ and $1\leq k\leq r,$%
\begin{eqnarray*}
\left\Vert W_{r}^{\ast }\left( X_{i}^{r}\right) ^{\infty }W_{r}e_{k}-\varphi
\left( x_{i}\right) e_{k}\right\Vert &<&\frac{1}{r}; \\
\text{ }\left\Vert W_{r}^{\ast }\left( Y_{i}^{r}\right) ^{\infty
}W_{r}e_{k}-\varphi \left( y_{i}\right) e_{k}\right\Vert &<&\frac{1}{r}\text{%
;} \\
\left\Vert W_{r}^{\ast }\left( Z_{j}^{r}\right) ^{\infty }W_{r}e_{k}-\varphi
\left( z_{j}\right) e_{k}\right\Vert &<&\frac{1}{r}.
\end{eqnarray*}

Therefore 
\begin{equation*}
W_{r}^{\ast }\left( X_{i}^{r}\right) ^{\left( \infty \right)
}W_{r}\rightarrow \varphi \left( x_{i}\right) \text{ \qquad in SOT as }%
r\rightarrow \infty
\end{equation*}%
\begin{equation*}
W_{r}^{\ast }\left( Y_{i}^{r}\right) ^{\left( \infty \right)
}W_{r}\rightarrow \varphi \left( y_{i}\right) \qquad \text{ in SOT as }%
r\rightarrow \infty
\end{equation*}%
\begin{equation*}
W_{r}^{\ast }\left( Z_{j}^{r}\right) ^{\left( \infty \right)
}W_{r}\rightarrow \varphi \left( z_{j}\right) \text{ \qquad in SOT as }%
r\rightarrow \infty
\end{equation*}%
for $i\in \mathbb{N},$ $j\in \left\{ 1,\cdots ,l\right\} $ and 
\begin{align*}
\left\Vert P(x_{1},x_{2},\cdots ,z_{1},\cdots ,z_{l})\right\Vert &
=\lim_{r\rightarrow \infty }\left\Vert P\left( X_{1}^{r},X_{2}^{r},\cdots
,Z_{1}^{r},\cdots ,Z_{l}^{r}\right) \right\Vert \\
\left\Vert Q(y_{1},y_{2},\cdots ,z_{1},\cdots ,z_{l})\right\Vert &
=\lim_{r\rightarrow \infty }\left\Vert Q\left( Y_{1}^{r},Y_{2}^{r},\cdots
,Z_{1}^{r},\cdots ,Z_{l}^{r}\right) \right\Vert
\end{align*}%
for any $P\in \mathbb{C}_{\mathbb{Q}}\left\langle \mathbf{X}_{1},\mathbf{X}%
_{2}\cdots ,\mathbf{Z}_{1},\cdots ,\mathbf{Z}_{l}\right\rangle $ and $Q\in 
\mathbb{C}_{\mathbb{Q}}\left\langle \mathbf{Y}_{1},\mathbf{Y}_{2}\cdots ,%
\mathbf{Z}_{1},\cdots ,\mathbf{Z}_{l}\right\rangle $ as desired.
\end{proof}

Now we are ready to show our main result in this subsection.

\begin{theorem}
\label{16}Let $\mathcal{A}$ and $\mathcal{B}$ be unital MF-algebras and $%
\mathcal{D}$ be a finite-dimensional C*-algebra. Suppose $\psi _{1}:\mathcal{%
D\rightarrow A}$ and $\psi _{2}:$ $\mathcal{D\rightarrow B}$ are unital
embeddings. Then $\mathcal{A}\underset{\mathcal{D}}{\mathcal{\ast }}\mathcal{%
B}$ is an MF algebra if and only if there are unital embeddings $q_{1}:%
\mathcal{A\rightarrow }\prod_{n=1}^{\infty }\mathcal{M}_{k_{n}}(\mathbb{C})/%
\overset{\infty }{\underset{n=1}{\sum }}\mathcal{M}_{k_{n}}(\mathbb{C})$ and 
$q_{2}:\mathcal{B\rightarrow }\prod_{n=1}^{\infty }\mathcal{M}_{k_{n}}(%
\mathbb{C})/\overset{\infty }{\underset{n=1}{\sum }}\mathcal{M}_{k_{n}}(%
\mathbb{C})$ for a sequence $\left\{ k_{n}\right\} _{n=1}^{\mathcal{1}}$ of
integers such that the following diagram commutes
\end{theorem}

\begin{equation*}
\begin{array}{ccc}
\mathcal{D} & \overset{\psi _{\mathcal{A}}}{\hookrightarrow } & \mathcal{A}
\\ 
^{\psi _{\mathcal{B}}}\downarrow &  & {\Huge \downarrow }^{q_{1}} \\ 
\mathcal{B} & \overset{q_{2}}{\hookrightarrow } & \prod_{m=1}^{\infty }%
\mathcal{M}_{k_{m}}(\mathbb{C)}/\tsum \mathcal{M}_{k_{m}}(\mathbb{C)}%
\end{array}%
\end{equation*}

\begin{proof}
If $\mathcal{A}\underset{\mathcal{D}}{\mathcal{\ast }}\mathcal{B}$ is an MF
algebra, then there is a unital embedding 
\begin{equation*}
\Phi :\mathcal{A}\underset{\mathcal{D}}{\mathcal{\ast }}\mathcal{%
B\rightarrow }\prod_{n=1}^{\infty }\mathcal{M}_{k_{n}}(\mathbb{C})/\overset{%
\infty }{\underset{n=1}{\sum }}\mathcal{M}_{k_{n}}(\mathbb{C})
\end{equation*}
for a sequence $\left\{ k_{n}\right\} _{n=1}^{\infty }$ of integers. Let $%
q_{1}$ and $q_{2}$ be the restrictions of $\Phi $ on $\mathcal{A}$ and $%
\mathcal{B}$ respectively. Then the above diagram is commutative.

Conversely, suppose $\mathcal{D},$ $\mathcal{A}$ and $\mathcal{B}$ are
generated by families%
\begin{equation*}
\left\{ z_{1},\cdots ,z_{l}\right\} ,
\end{equation*}%
\begin{equation*}
\left\{ x_{1},x_{2}\cdots ,\psi _{1}\left( z_{1}\right) ,\cdots ,\psi
_{1}\left( z_{l}\right) \right\}
\end{equation*}%
and 
\begin{equation*}
\left\{ y_{1},y_{2},\cdots ,\psi _{2}\left( z_{1}\right) ,\cdots ,\psi
_{2}\left( z_{l}\right) \right\}
\end{equation*}%
respectively with $\left\Vert x_{i}\right\Vert =\left\Vert y_{i}\right\Vert
=\left\Vert z_{j}\right\Vert =1$ for each $i\in \mathbb{N}$ and $j\in
\left\{ 1,\cdots ,l\right\} $. By Remark \ref{13.3}, we may assume that $%
\mathcal{D}$ is a finite-dimensional abelian C*-algebra and $z_{1},\cdots
z_{l}$ are orthogonal projections with $\tsum\limits_{i=1}^{l}z_{i}=I$.
Without loss of generality, we may assume that $\mathcal{A}\underset{%
\mathcal{D}}{\mathcal{\ast }}\mathcal{B}$ is generated by a sequence $%
\left\{ x_{1},x_{2}\cdots ,y_{1},y_{2}\cdots ,z_{1},\cdots ,z_{l}\right\} .$
Let $\varphi :\mathcal{A}\underset{\mathcal{D}}{\mathcal{\ast }}\mathcal{%
B\rightarrow B}\left( \mathcal{H}\right) $ be a faithful representation of
full amalgamated free product $\mathcal{A}\underset{\mathcal{D}}{\mathcal{%
\ast }}\mathcal{B}$. Applying Proposition \ref{15}, there is a sequence $%
\left\{ t_{m}\right\} _{m=1}^{\mathcal{1}}$ of integers such that, for each $%
t_{r}\in \left\{ t_{m}\right\} _{m=1}^{\mathcal{1}},$ there exist sequences $%
\left\{ X_{1}^{r},X_{2}^{r},\cdots \right\} ,\left\{
Y_{1}^{r},Y_{2}^{r},\cdots \right\} $ and $\left\{ Z_{1}^{r},\cdots
,Z_{2}^{r}\right\} $ in $\mathcal{M}_{t_{r}}\left( \mathbb{C}\right) $ and a
unitary operator $W_{r}:\mathcal{H\rightarrow }\left( \mathbb{C}%
^{t_{r}}\right) ^{\infty }$ such that 
\begin{equation*}
W_{r}^{\ast }\left( X_{i}^{r}\right) ^{\left( \infty \right)
}W_{r}\rightarrow \varphi \left( x_{i}\right) \text{ \qquad in SOT as }%
r\rightarrow \infty \text{ for }i\in \mathbb{N}
\end{equation*}%
\begin{equation*}
W_{r}^{\ast }\left( Y_{i}^{r}\right) ^{\left( \infty \right)
}W_{r}\rightarrow \varphi \left( y_{i}\right) \text{ \qquad in SOT as }%
r\rightarrow \infty \text{ for }i\in \mathbb{N}
\end{equation*}%
and%
\begin{equation*}
W_{r}^{\ast }\left( Z_{j}^{r}\right) ^{\left( \infty \right)
}W_{r}\rightarrow \varphi \left( z_{j}\right) \qquad \text{ in SOT as }%
r\rightarrow \infty \text{ for }j\in \left\{ 1,\cdots ,l\right\}
\end{equation*}%
as well as 
\begin{align*}
\left\Vert P\left( x_{1},x_{2},\cdots ,z_{1},z_{2},\cdots \right)
\right\Vert & =\lim_{r\rightarrow \infty }\left\Vert P\left(
X_{1}^{r},X_{2}^{r},\cdots ,Z_{1}^{r},\cdots ,Z_{l}^{r}\right) \right\Vert \\
\left\Vert Q\left( y_{1},y_{2},\cdots ,z_{1},z_{2},\cdots \right)
\right\Vert & =\lim_{r\rightarrow \infty }\left\Vert Q\left(
Y_{1}^{r},Y_{2}^{r},\cdots ,Z_{1}^{r},\cdots ,Z_{l}^{r}\right) \right\Vert
\end{align*}%
for any $P\in \mathbb{C}_{\mathbb{Q}}\left\langle \mathbf{X}_{1},\mathbf{X}%
_{2}\cdots ,\mathbf{Z}_{1},\cdots ,\mathbf{Z}_{l}\right\rangle $ and $Q\in 
\mathbb{C}_{\mathbb{Q}}\left\langle \mathbf{Y}_{1},\mathbf{Y}_{2}\cdots ,%
\mathbf{Z}_{1},\cdots ,\mathbf{Z}_{l}\right\rangle .$ Therefore, we can
define unital embeddings 
\begin{equation*}
q_{1}:\mathcal{A}\rightarrow \tprod M_{t_{r}}\left( \mathbb{C}\right) /\tsum
M_{t_{r}}\left( \mathbb{C}\right)
\end{equation*}%
and 
\begin{equation*}
q_{2}:\mathcal{B}\rightarrow \tprod M_{t_{r}}\left( \mathbb{C}\right) /\tsum
M_{t_{r}}\left( \mathbb{C}\right)
\end{equation*}%
so that $q_{1}(x_{i})=\left[ \left( X_{i}^{r}\right) \right] ,q_{1}(\psi
_{1}\left( z_{j}\right) )=\left[ \left( Z_{j}^{r}\right) \right] ,$ and $%
q_{2}(y_{i})=\left[ \left( Y_{i}^{r}\right) \right] ,q_{2}(\psi _{2}\left(
z_{j}\right) )=\left[ \left( Z_{j}^{r}\right) \right] $ for $i\in \mathbb{N}%
,j\in \left\{ 1,\cdots ,l\right\} .$ From the definition of full amalgamated
free product, there is a $\ast $-homomorphism 
\begin{equation*}
\Phi :\mathcal{A}\underset{\mathcal{D}}{\mathcal{\ast }}\mathcal{%
B\rightarrow }\tprod M_{t_{r}}\left( \mathbb{C}\right) /\tsum
M_{t_{r}}\left( \mathbb{C}\right)
\end{equation*}%
such that $\Phi (x_{i})=\left[ \left( X_{i}^{r}\right) \right] ,\Phi (y_{i})=%
\left[ \left( Y_{i}^{r}\right) \right] ,\Phi (z_{j})=\left[ \left(
Z_{j}^{r}\right) \right] $ where $i\in \mathbb{N}$ and $j\in \left\{
1,\cdots ,l\right\} .$ Furthermore, for any $\Psi _{j}\in \mathbb{C}_{%
\mathbb{Q}}\left\langle \mathbf{X}_{1},\mathbf{X}_{2}\cdots ,\mathbf{Y}_{1},%
\mathbf{Y}_{2}\cdots ,\mathbf{Z}_{1},\cdots ,\mathbf{Z}_{l}\right\rangle ,$
we have%
\begin{align}
& \left\Vert \Psi _{j}\left( \left[ \left( X_{1}^{r}\right) \right] ,\left[
\left( X_{2}^{r}\right) \right] \cdots ,\left[ \left( Y_{1}^{r}\right) %
\right] ,\left[ \left( Y_{2}^{r}\right) \right] \cdots ,\left[ \left(
Z_{1}^{r}\right) \right] ,\cdots ,\left[ \left( Z_{l}^{r}\right) \right]
\right) \right\Vert  \notag \\
& =\underset{r\rightarrow \infty }{\lim \sup }\left\Vert \Psi _{j}\left(
X_{1}^{r},X_{2}^{r}\cdots ,Y_{1}^{r},Y_{2}^{r}\cdots ,,Z_{1}^{r},\cdots
,Z_{l}^{r}\right) \right\Vert _{M_{k_{r}}\left( \mathbb{C}\right) }  \notag
\\
& \leq \left\Vert \Psi _{j}\left( x_{1},x_{2}\cdots ,y_{1},y_{2}\cdots
,z_{1},\cdots ,z_{l}\right) \right\Vert _{\mathcal{A}\underset{\mathcal{D}}{%
\mathcal{\ast }}\mathcal{B}}.\qquad \ \ \ \ \ \ \ \ \   \tag{5.38}
\label{38}
\end{align}%
Meanwhile, 
\begin{align*}
& \Psi _{j}\left( W_{r}^{\ast }\left( X_{1}^{r}\right) ^{\infty
}W_{r},\cdots ,W_{r}^{\ast }\left( Y_{1}^{r}\right) ^{\infty }W_{r},\cdots
,W_{r}^{\ast }\left( Z_{1}^{r}\right) ^{\mathcal{1}}W_{r},\cdots
,W_{r}^{\ast }\left( Z_{l}^{r}\right) ^{\mathcal{1}}W_{r}\right) \\
& \rightarrow \Psi _{j}\left( x_{1},x_{2}\cdots ,y_{1},y_{2}\cdots
,z_{1},\cdots ,z_{l}\right) \text{ in SOT as }r\rightarrow \infty ,
\end{align*}%
and therefore%
\begin{align}
& \underset{r\rightarrow \infty }{\lim \inf }\left\Vert \Psi _{j}\left(
X_{1}^{r},X_{2}^{r}\cdots ,Y_{1}^{r},Y_{2}^{r}\cdots ,,Z_{1}^{r},\cdots
,Z_{l}^{r}\right) \right\Vert _{M_{k_{r}}\left( \mathbb{C}\right) }  \notag
\\
& =\underset{r\rightarrow \infty }{\lim \inf }\left\Vert \Psi _{j}\left(
W_{r}^{\ast }\left( X_{1}^{r}\right) ^{\infty }W_{r},\cdots ,W_{r}^{\ast
}\left( Y_{1}^{r}\right) ^{\infty }W_{r},\cdots ,W_{r}^{\ast }\left(
Z_{1}^{r}\right) ^{\mathcal{1}}W_{r},\cdots ,W_{r}^{\ast }\left(
Z_{l}^{r}\right) ^{\mathcal{1}}W_{r}\right) \right\Vert  \notag \\
& \geq \left\Vert \Psi _{j}\left( x_{1},x_{2}\cdots ,y_{1},y_{2}\cdots
,z_{1},\cdots ,z_{l}\right) \right\Vert _{\mathcal{A}\underset{\mathcal{D}}{%
\mathcal{\ast }}\mathcal{B}}.\qquad \qquad \qquad \qquad \qquad \ \ \ \ \ \
\ \ \ \ \ \ \ \ \ \ \text{\ }\   \tag{5.39}  \label{39}
\end{align}%
Combining (\ref{38}) and (\ref{39})$,$ it follows that 
\begin{equation*}
\left\Vert \Psi \left( x_{1},x_{2}\cdots ,y_{1},y_{2}\cdots ,z_{1},\cdots
,z_{l}\right) \right\Vert =\underset{r\rightarrow \infty }{\lim }\left\Vert
\Psi _{j}\left( X_{1}^{r},X_{2}^{r}\cdots ,Y_{1}^{r},Y_{2}^{r}\cdots
,Z_{1}^{r},\cdots ,Z_{l}^{r}\right) \right\Vert
\end{equation*}%
for any $\Psi \in \mathbb{C}_{\mathbb{Q}}\left\langle \mathbf{X}_{1},\cdots ,%
\mathbf{Y}_{1},\cdots ,\mathbf{Z}_{1},\cdots ,\mathbf{Z}_{l}\right\rangle .$
Then $\Phi $ is a unital injective $\ast $- homomorphism. It follows that $%
\mathcal{A}\underset{\mathcal{D}}{\mathcal{\ast }}\mathcal{B}$ is an MF
algebra.
\end{proof}

The following corollary is an easy consequence of Theorem \ref{16}.

\begin{corollary}
\label{17}Let $\mathcal{A}$ be an MF algebra and $\mathcal{D}$ be a
finite-dimensional C*-algebra. If there is a unital embedding $q:\mathcal{%
D\rightarrow A},$ then $\mathcal{A}\underset{\mathcal{D}}{\mathcal{\ast }}%
\mathcal{A}$ is an MF algebra with respect to the embedding $q.$
\end{corollary}

Applying Theorems \ref{16}, we can obtain the following result.

\begin{proposition}
\label{19}Let $\mathcal{A}$ and $\mathcal{B}$ be unital MF algebras and let $%
\mathcal{D}$ be the direct sum of n copies of the set of all complex
numbers, that is, 
\begin{equation*}
\mathcal{D}=\mathbb{C}\oplus \mathbb{C}\oplus \cdots \oplus \mathbb{C}.
\end{equation*}%
Suppose $\psi _{\mathcal{A}}:\mathcal{D\rightarrow A}$ and $\psi _{\mathcal{B%
}}:\mathcal{D\rightarrow B}$ are unital embeddings. If $\psi _{\mathcal{A}}$
and $\psi _{\mathcal{B}}$ can be extended to unital embeddings $\widetilde{%
\psi }_{\mathcal{A}}:$ $\mathcal{M}_{n}\left( \mathbb{C}\right) \rightarrow 
\mathcal{A}$ and $\widetilde{\psi }_{\mathcal{B}}:$ $\mathcal{M}_{n}\left( 
\mathbb{C}\right) \rightarrow \mathcal{B}$ respectively, then $\mathcal{A}%
\underset{\mathcal{D}}{\mathcal{\ast }}\mathcal{B}$ is an MF algebra.\textbf{%
\ }
\end{proposition}

\begin{proof}
Let 
\begin{equation*}
E_{1}=1\oplus 0\oplus \cdots \oplus 0,E_{2}=0\oplus 1\oplus 0\oplus \cdots
\oplus 0,\cdots ,E_{n}=0\oplus 0\oplus \cdots \oplus 0\oplus 1
\end{equation*}%
in $\mathcal{D}.$ Then $\mathcal{D=}\mathbb{C}E_{1}+\cdots +\mathbb{C}E_{n}.$
Suppose C*-algebras $\mathcal{A}$ and $\mathcal{B}$ are generated by
families 
\begin{equation*}
\left\{ x_{1},x_{2}\cdots ,\psi _{\mathcal{A}}\left( E_{1}\right) ,\cdots
,\psi _{\mathcal{A}}\left( E_{n}\right) \right\}
\end{equation*}%
of self-adjoint elements and 
\begin{equation*}
\left\{ y_{1},y_{2}\cdots ,\psi _{\mathcal{B}}\left( E_{1}\right) ,\cdots
,\psi _{\mathcal{B}}\left( E_{n}\right) \right\}
\end{equation*}%
of self-adjoint elements respectively. Without loss of generality, we may
assume that $\mathcal{A}$ and $\mathcal{B}$ can be embedded as unital
C*-subalgebras of $\prod_{m=1}^{\infty }\mathcal{M}_{k_{m}}(\mathbb{C)}%
/\tsum \mathcal{M}_{k_{m}}(\mathbb{C)}$, respectively, for a sequence $%
\left\{ k_{m}\right\} _{m=1}^{m}$ of integers with sequences 
\begin{equation*}
\left\{ A_{1}^{m},A_{2}^{m},\cdots \right\} ,\left\{ C_{1}^{m},\cdots
,C_{n}^{m}\right\} \text{ and }\left\{ B_{1}^{m},B_{2}^{m},\cdots \right\}
,\left\{ D_{1}^{m},\cdots ,D_{n}^{m}\right\} \subset \mathcal{M}%
_{k_{m}}^{s.a.}\left( \mathbb{C}\right)
\end{equation*}%
for each $k_{m}\in \left\{ k_{n}\right\} _{n=1}^{\mathcal{1}}$ satisfying 
\begin{equation*}
\lim_{m\rightarrow \infty }\Vert P(A_{1}^{m},A_{2}^{m}\ldots
,C_{1}^{m},\cdots ,C_{n}^{m})\Vert =\Vert P(x_{1},x_{2}\cdots ,\psi _{%
\mathcal{A}}\left( E_{1}\right) ,\cdots ,\psi _{\mathcal{A}}\left(
E_{n}\right) )\Vert
\end{equation*}%
for any $P\in \mathbb{C}_{\mathbb{Q}}\left\langle \mathbf{X}_{1},\mathbf{X}%
_{2}\cdots ,\mathbf{Z}_{1},\cdots ,\mathbf{Z}_{n}\right\rangle $, and 
\begin{equation*}
\lim_{m\rightarrow \infty }\Vert Q(B_{1}^{m},B_{2}^{m}\ldots
,D_{1}^{m},\cdots ,D_{n}^{m})\Vert =\Vert Q(y_{1},y_{2}\cdots ,\psi _{%
\mathcal{B}}\left( E_{1}\right) ,\cdots ,\psi _{\mathcal{B}}\left(
E_{n}\right) )\Vert
\end{equation*}%
for any $Q\in \mathbb{C}_{\mathbb{Q}}\left\langle \mathbf{Y}_{1},\mathbf{Y}%
_{2},\cdots ,\mathbf{Z}_{1},\cdots ,\mathbf{Z}_{n}\right\rangle .$ Since the
images of $\psi _{\mathcal{A}}\left( E_{1}\right) ,\cdots ,\psi _{\mathcal{A}%
}\left( E_{n}\right) $ under the embedding from $\mathcal{A}$ to $%
\prod_{m=1}^{\infty }\mathcal{M}_{k_{m}}(\mathbb{C)}/\tsum \mathcal{M}%
_{k_{m}}(\mathbb{C)}$ are $\left[ \left( C_{1}^{m}\right) \right] ,\cdots ,%
\left[ \left( C_{n}^{m}\right) \right] $ and $\psi _{\mathcal{A}}$ can be
extended to a unital embedding $\widetilde{\psi }_{\mathcal{A}}:$ $\mathcal{M%
}_{n}\left( \mathbb{C}\right) \rightarrow \mathcal{A}$, it follows that
there are partial isometries $\left[ \left( V_{1}^{m}\right) \right] ,\cdots
,\left[ \left( V_{n}^{m}\right) \right] $ in $\prod_{m=1}^{\infty }\mathcal{M%
}_{k_{m}}(\mathbb{C)}/\tsum \mathcal{M}_{k_{m}}(\mathbb{C)}$ such that $%
\left[ \left( V_{s}^{m}\right) \right] ^{\ast }\left[ \left(
V_{s}^{m}\right) \right] =\left[ \left( C_{1}^{m}\right) \right] $ and $%
\left[ \left( V_{s}^{m}\right) \right] \left[ \left( V_{s}^{m}\right) \right]
^{\ast }=\left[ \left( C_{s}^{m}\right) \right] $ for each $s\in \left\{
1,\cdots ,n\right\} .$ By Lemma \ref{18}, we may assume that, $C_{s}^{m}\in 
\mathcal{M}_{k_{m}}(\mathbb{C)}$ is a projection for each $m\in \mathbb{N}$
and $s\in \left\{ 1,\cdots ,n\right\} .$ We may conclude further that, when $%
k_{m}$ is large enough, $V_{s}^{m}$ is a partial isometry such that $%
V_{s}^{m\ast }V_{s}^{m}=C_{1}^{m}$ and $V_{s}^{m}V_{s}^{m\ast }=C_{s}^{m}$
for $1\leq s\leq n$ in $\mathcal{M}_{k_{m}}\left( \mathbb{C}\right) $ by
Lemma \ref{18}. So it follows that $C_{1}^{m}$ is equivalent to $C_{s}^{m}$
in $\mathcal{M}_{k_{m}}(\mathbb{C)}$ for $1\leq s\leq n$ and $%
\tsum\limits_{s=1}^{n}C_{s}^{m}=I$, $C_{i}^{m}C_{j}^{m}=0$ for $1\leq i\neq
j\leq n.$ Similarly, we can assume that $D_{s}^{m}$ is a projection in $%
\mathcal{M}_{k_{m}}(\mathbb{C)}$ for $1\leq s\leq n$. When $k_{m}$ is large
enough, we conclude that $D_{1}^{m}$ is equivalent to $D_{s}^{m}$ for each $%
1\leq s\leq n,$ $\tsum_{s=1}^{n}$ $D_{s}^{m}=I$ and $D_{i}^{m}D_{j}^{m}=0$
as $1\leq i\neq j\leq n$ in $\mathcal{M}_{k_{m}}(\mathbb{C)}.$ Hence, there
exists an integer $K$ such that, for each $k_{m}>K,$ there exists a unitary $%
U^{m}\in \mathcal{M}_{k_{m}}\left( \mathbb{C}\right) $ satisfying $%
U^{m}C_{s}^{m}U^{m\ast }=D_{s}^{m}$ for each $s\in \left\{ 1,2,\cdots
,n\right\} $ in $\mathcal{M}_{k_{m}}(\mathbb{C)}$. It follows that there is
a unitary $\left[ \left( U^{m}\right) \right] \in \prod_{m=1}^{\infty }%
\mathcal{M}_{k_{m}}(\mathbb{C)}/\tsum \mathcal{M}_{k_{m}}(\mathbb{C)}$
satisfying $\left[ \left( U^{m}\right) \right] \left[ \left(
C_{i}^{m}\right) \right] \left[ \left( U^{m}\right) \right] ^{\ast }=\left[
D_{i}^{m}\right] $ for $1\leq i\leq n.$ Now we define embeddings 
\begin{equation*}
q_{1}:\mathcal{A\rightarrow }\prod_{m=1}^{\infty }\mathcal{M}_{k_{m}}(%
\mathbb{C)}/\tsum \mathcal{M}_{k_{m}}(\mathbb{C)}
\end{equation*}%
so that $q_{1}\left( x_{i}\right) =\left[ \left( U^{m}\right) \right] \left[
\left( A_{i}^{m}\right) \right] \left[ \left( U^{m}\right) \right] ^{\ast }$
for $i\in \mathbb{N}$, $q_{1}\left( \psi _{\mathcal{A}}\left( z_{j}\right)
\right) =\left[ \left( D_{j}^{m}\right) \right] $ for $1\leq j\leq n$ and 
\begin{equation*}
q_{2}:\mathcal{B\rightarrow }\prod_{m=1}^{\infty }\mathcal{M}_{k_{m}}(%
\mathbb{C)}/\tsum \mathcal{M}_{k_{m}}(\mathbb{C)}
\end{equation*}%
so that $q_{2}\left( y_{i}\right) =\left[ \left( B_{i}^{m}\right) \right] $
for $i\in \mathbb{N}$, $q_{2}\left( \psi _{\mathcal{B}}\left( z_{i}\right)
\right) =\left[ \left( D_{j}^{m}\right) \right] $ for$1\leq j\leq n.$ It is
clear that the following diagram is commutative%
\begin{equation*}
\begin{array}{ccc}
\mathcal{D} & \overset{\psi _{\mathcal{A}}}{\hookrightarrow } & \mathcal{A}
\\ 
^{\psi _{\mathcal{B}}}\downarrow &  & {\Huge \downarrow }^{q_{1}} \\ 
\mathcal{B} & \overset{q_{2}}{\hookrightarrow } & \prod_{m=1}^{\infty }%
\mathcal{M}_{k_{m}}(\mathbb{C)}/\tsum \mathcal{M}_{k_{m}}(\mathbb{C)}%
\end{array}%
\end{equation*}%
So $\mathcal{A}\underset{\mathcal{D}}{\mathcal{\ast }}\mathcal{B}$ is MF by
Theorem \ref{16}.
\end{proof}

\subsection{$\mathcal{D}$ Is An Infinite-dimensional C*-algebra}

In this subsection, we consider the case when $\mathcal{D}$ is an
infinite-dimensional C*-algebra. More precisely, we consider the case when $%
\mathcal{D}$ can be written as a norm closure of the union of an increasing
sequence of C*-algebras.

\begin{theorem}
\label{Theorem 3.1} Suppose that $\mathcal{A}\supseteq \mathcal{D}\subseteq 
\mathcal{B}$ are unital inclusions of unital separable C*-algebras and $\{%
\mathcal{D}_{k}\}_{k=1}^{\infty }$ is an increasing sequence of unital C$%
^{\ast }$-subalgebras of $\mathcal{D}$ such that $\cup _{k\geq 1}\mathcal{D}%
_{k}$ is norm dense in $\mathcal{D}$. Let $\mathcal{A}\ast _{\mathcal{D}}%
\mathcal{B}$ and $\mathcal{A}\ast _{\mathcal{D}_{k}}\mathcal{B}$ for $k\geq
1 $ be the unital full free products of $\mathcal{A}$ and $\mathcal{B}$ with
amalgamation over $\mathcal{D\ }$and $\mathcal{D}_{k}$ for $k\geq 1$
respectively. If $\mathcal{A}\ast _{\mathcal{D}_{k}}\mathcal{B}$ is an MF
algebra for each $k\geq 1$, then $\mathcal{A}\ast _{\mathcal{D}}\mathcal{B}$
is an MF algebra.
\end{theorem}

\begin{proof}
Note that $\mathcal{A},\mathcal{B}$ and $\mathcal{D}$ are unital separable C$%
^{\ast }$-algebras. We might assume that $\{x_{n}\}_{n\in \mathbb{N}%
}\subseteq \mathcal{A}$, and $\{y_{n}\}_{n\in \mathbb{N}}\subseteq \mathcal{B%
}$ are families of generators of $\mathcal{A}$ and $\mathcal{B}$
respectively.

Assume that $\sigma :\mathcal{A}\rightarrow \mathcal{A}\ast _{\mathcal{D}}%
\mathcal{B}$ and $\sigma _{k}:\mathcal{A}\rightarrow \mathcal{A}\ast _{%
\mathcal{D}_{k}}\mathcal{B}$ are natural unital embeddings from $\mathcal{A}$
into $\mathcal{A}\ast _{\mathcal{D}}\mathcal{B}$ and $\mathcal{A}\ast _{%
\mathcal{D}_{k}}\mathcal{B}$, respectively, for each $k\geq 1$. Assume that $%
\rho :\mathcal{B}\rightarrow \mathcal{A}\ast _{\mathcal{D}}\mathcal{B}$ and $%
\rho _{k}:\mathcal{B}\rightarrow \mathcal{A}\ast _{\mathcal{D}_{k}}\mathcal{B%
}$ are natural unital embeddings from $\mathcal{B}$ into $\mathcal{A}\ast _{%
\mathcal{D}}\mathcal{B}$ and into $\mathcal{A}\ast _{\mathcal{D}_{k}}%
\mathcal{B}$, respectively, for each $k\geq 1$.

Consider the unital C$^{\ast }$-algebra 
\begin{equation*}
\prod_{k\geq 1}\mathcal{A}\ast _{\mathcal{D}_{k}}\mathcal{B}/\sum_{k\geq 1}%
\mathcal{A}\ast _{\mathcal{D}_{k}}\mathcal{B}.
\end{equation*}%
From Corollary 3.4.3 in \cite{[BK]} and the fact that, for each $k\geq 1$, $%
\mathcal{A}\ast _{\mathcal{D}_{k}}\mathcal{B}$ is an MF algebra, we know
that every separable C$^{\ast }$-subalgebra of $\prod_{k\geq 1}\mathcal{A}%
\ast _{\mathcal{D}}\mathcal{B}/\prod_{k\geq 1}\mathcal{A}\ast _{\mathcal{D}%
_{k}}\mathcal{B}$ is an MF algebra. Let 
\begin{equation*}
a_{n}=[(\sigma _{k}(x_{n}))_{k}]\in \prod_{k\geq 1}\mathcal{A}\ast _{%
\mathcal{D}_{k}}\mathcal{B}/\sum_{k\geq 1}\mathcal{A}\ast _{\mathcal{D}_{k}}%
\mathcal{B},\qquad \forall \ n\in \mathbb{N}
\end{equation*}%
and 
\begin{equation*}
b_{n}=[(\rho _{k}(y_{n}))_{k}]\in \prod_{k\geq 1}\mathcal{A}\ast _{\mathcal{D%
}_{k}}\mathcal{B}/\sum_{k\geq 1}\mathcal{A}\ast _{\mathcal{D}_{k}}\mathcal{B}%
,\qquad \forall \ n\in \mathbb{N}.
\end{equation*}%
Let $\mathcal{Q}$ be the unital C$^{\ast }$-subalgebra generated by $%
\{a_{n},b_{n}\}_{n\in \mathbb{N}}$ in $\prod_{k\geq 1}\mathcal{A}\ast _{%
\mathcal{D}_{k}}\mathcal{B}/\sum_{k\geq 1}\mathcal{A}\ast _{\mathcal{D}_{k}}%
\mathcal{B}$. Thus $\mathcal{Q}$ is an MF algebra.

Next we shall show that there is a $\ast $-isomorphic from $\mathcal{Q}$
onto $\mathcal{A}\ast _{\mathcal{D}}\mathcal{B}$ by sending each $a_{n}$ to $%
\sigma (x_{n})$ and $b_{n}$ to $\rho (y_{n})$. This will induce that $%
\mathcal{A}\ast _{\mathcal{D}}\mathcal{B}$ is also an MF algebra. In order
to obtain such $\ast $-isomorphism from $\mathcal{Q}$ onto $\mathcal{A}\ast
_{\mathcal{D}}\mathcal{B}$, it suffices to show that \emph{$\forall \ N\in 
\mathbb{N}$ and $\forall \ P\in \mathbb{C\langle }X_{1},\mathbb{\ldots }%
.X_{N},Y_{1},\mathbb{\ldots }Y_{N}\mathbb{\rangle }$, we have 
\begin{eqnarray}
&&\left\Vert P(\sigma (x_{1}),\cdots ,\sigma (x_{N}),\rho (y_{1}),\cdots
,\rho (y_{N}))\right\Vert _{\mathcal{A\ast }_{\mathcal{D}}\mathcal{B}}\ \ \
\ \ \ \ \ \ \ \ \ \ \ \ \ \ \ \ \ \ \ \ \ \   \notag \\
\ &=&\left\Vert P(a_{1},\cdots ,a_{N},b_{1},\cdots ,b_{N})\right\Vert
_{\tprod\limits_{k\geq 1}\mathcal{A\ast }_{\mathcal{D}_{k}}\mathcal{B}/\sum 
\mathcal{A\ast }_{\mathcal{D}_{k}}\mathcal{B}}  \TCItag{$5.40$}
\label{eq3.1}
\end{eqnarray}%
}

By the definition of full amalgamated free product, we know, for each $k\geq
1,$ there is a $\ast $-homomorphism from $\mathcal{A}\ast _{\mathcal{D}_{k}}%
\mathcal{B}$ to $\mathcal{A}\ast _{\mathcal{D}}\mathcal{B}$, which send $%
\sigma _{k}(x_{n})\ $to $\sigma (x_{n})$ and $\rho _{k}(y_{n})$ to $\rho
(y_{n})$ respectively, for every $n\in \mathbb{N}$. Hence 
\begin{eqnarray*}
&&\left\Vert P(\sigma (x_{1}),\cdots ,\sigma (x_{N}),\rho (y_{1}),\cdots
,\rho (y_{N}))\right\Vert _{\mathcal{A\ast }_{\mathcal{D}}\mathcal{B}} \\
&\leq &\left\Vert P(\sigma _{k}(x_{1}),\cdots ,\sigma _{k}(x_{N}),\rho
_{k}(y_{1}),\cdots ,\rho _{k}(y_{N}))\right\Vert _{\mathcal{A\ast }_{%
\mathcal{D}_{k}}\mathcal{B}}\ \ \text{for\ all\ }k\geq 1
\end{eqnarray*}%
and, consequently, 
\begin{equation}
\left\Vert P(\sigma (x_{1}),\cdots ,\sigma (x_{N}),\rho (y_{1}),\cdots ,\rho
(y_{N}))\right\Vert _{\mathcal{A\ast }_{\mathcal{D}}\mathcal{B}}\leq
\left\Vert P(a_{1},\cdots ,a_{N},b_{1},\cdots ,b_{N})\right\Vert
_{\tprod\limits_{k\geq 1}\mathcal{A}\ast _{\mathcal{D}_{k}}\mathcal{B}%
/\sum_{k\geq 1}\mathcal{A\ast }_{\mathcal{D}_{k}}\mathcal{B}}  \tag{5.41}
\label{eq3.2}
\end{equation}

We will show that $[(\sigma _{k}(z))_{k}]=[(\rho _{k}(z))_{k}]$ for every $%
z\in \mathcal{D}$. Suppose $z\in \mathcal{D}$ and $\varepsilon >0$. Then
there exist a positive integer $p$ and an element $z_{p}$ in $\mathcal{D}%
_{p} $ such that $\Vert z-z_{p}\Vert <\varepsilon .$ Since $\{\mathcal{D}%
_{k}\}_{k}$ is an increasing sequence of C$^{\ast }$-algebras, we know that $%
z_{p}\in \mathcal{D}_{k}$ for $k\geq p$. It follows that $[(\sigma
_{k}(z_{p}))_{k}]=[(\rho _{k}(z_{p}))_{k}]$. So,%
\begin{eqnarray*}
&&\left\Vert [(\sigma _{k}\left( z\right) )_{k}]-[(\rho _{k}\left( z\right)
)_{k}]\right\Vert _{\tprod\limits_{k\geq 1}\mathcal{A\ast }_{\mathcal{D}_{k}}%
\mathcal{B}/\sum_{k\geq 1}\mathcal{A\ast }_{\mathcal{D}_{k}}\mathcal{B}} \\
&\mathcal{=}&\underset{k}{\lim \sup }\left\Vert \sigma _{k}\left( z\right)
-\rho _{k}\left( z\right) \right\Vert _{\mathcal{A\ast }_{\mathcal{D}_{k}}}
\\
&\mathcal{\leq }&\underset{k}{\lim \sup }(\left\Vert \sigma _{k}(z)-\sigma
_{k}(z_{p})\right\Vert _{\mathcal{A\ast }_{\mathcal{D}_{k}}\mathcal{B}%
}+\left\Vert \sigma _{k}(z_{p})-\rho _{k}(z_{p})\right\Vert _{\mathcal{A\ast 
}_{\mathcal{D}_{k}}\mathcal{B}}+\left\Vert \rho _{k}(z_{p})-\rho
_{k}(z)\right\Vert _{\mathcal{A\ast }_{\mathcal{D}_{k}}\mathcal{B}}) \\
&\leq &2\varepsilon \ \ \ \ \ \ \text{for\ all\ }\varepsilon >0.
\end{eqnarray*}

Thus we obtain that $[(\sigma _{k}(z))_{k}]=[(\rho _{k}(z))_{k}]$.

Now it follows from the definitions of full amalgamated free product and of
the C$^{\ast }$-algebra $\mathcal{Q}$, together with the fact that $[(\sigma
_{k}(z))_{k}]=[(\rho _{k}(z))_{k}]$ for every $z\in \mathcal{D}$, we know
there is a $\ast $-homomorphism from $\mathcal{A}\ast _{\mathcal{D}}\mathcal{%
B}$ onto $\mathcal{Q}$, which maps each $\sigma (x_{n}),\rho (y_{n})$ to $%
a_{n},b_{n}$ respectively for $n\in \mathbb{N}$. Therefore,%
\begin{eqnarray}
&&\left\Vert P(\sigma (x_{1}),\cdots ,\sigma (x_{N}),\rho (y_{1}),\cdots
,\rho (y_{N}))\right\Vert _{\mathcal{A\ast }_{\mathcal{D}}\mathcal{B}} 
\notag \\
&\geq &\left\Vert P(a_{1},\cdots ,a_{N},b_{1},\cdots ,b_{N})\right\Vert
_{\tprod\limits_{k\geq 1}\mathcal{A\ast }_{\mathcal{D}_{k}}\mathcal{B}%
/\sum_{k\geq 1}\mathcal{A\ast }_{\mathcal{D}_{k}}\mathcal{B}}  \TCItag{5.42}
\label{eq3.3}
\end{eqnarray}

Now equation (\ref{eq3.1}) follows easily from inequalities (\ref{eq3.2})
and (\ref{eq3.3}). This ends our proof.
\end{proof}

Once we get the preceding theorem, we are ready to consider the case when $%
\mathcal{D}$ is an AF algebra. The following theorem is an analogous result
to Theorem \ref{16}

\begin{theorem}
\label{Theorem 3.2} Suppose that $\mathcal{A}\supset \mathcal{D}\subset 
\mathcal{B}$ are unital inclusions of unital MF algebras where $\mathcal{D}$
is an AF algebra. Then the unital full free product $\mathcal{A}\ast _{%
\mathcal{D}}\mathcal{B}$ of $\mathcal{A}$ and $\mathcal{B}$ with
amalgamation over $\mathcal{D}$ is an MF algebra if and only if there is a
sequence of positive integers $\{n_{k}\}_{k=1}^{\infty }$ such that the
following diagram 
\begin{equation*}
\begin{array}{ccc}
\mathcal{D} & \subseteq & \mathcal{A} \\ 
\cap &  & \cap \\ 
\mathcal{B} & \subseteq & \tprod\limits_{k}\mathcal{M}_{n_{k}}\left( \mathbb{%
C}\right) /\sum_{k}\mathcal{M}_{n_{k}}\left( \mathbb{C}\right)%
\end{array}%
\end{equation*}%
commutes. .
\end{theorem}

\begin{proof}
If $\mathcal{A}\ast _{\mathcal{D}}\mathcal{B}$ is MF, there is a sequence of
positive integers $\{n_{k}\}_{k=1}^{\infty }$ such that the following
diagram commutes automatically. For another direction, suppose there is a
sequence of positive integers $\{n_{k}\}_{k=1}^{\infty }$ such that the
following diagram 
\begin{equation*}
\begin{array}{ccc}
\mathcal{D} & \subseteq & \mathcal{A} \\ 
\cap &  & \cap \\ 
\mathcal{B} & \subseteq & \tprod\limits_{k}\mathcal{M}_{n_{k}}\left( \mathbb{%
C}\right) /\sum_{k}\mathcal{M}_{n_{k}}\left( \mathbb{C}\right)%
\end{array}%
\end{equation*}%
commutes. Note that $\mathcal{D}$ is an AF algebra, therefore there is an
increasing sequence of unital finite-dimensional C$^{\ast }$-subalgebra $\{%
\mathcal{D}_{p}\}_{p\geq 1}$ of $\mathcal{D}$ such that $\cup _{p}\mathcal{D}%
_{p}$ is norm dense in $\mathcal{D}$. By the choice of each $\mathcal{D}_{p}$%
, we know that the diagram, for each $p\geq 1$, 
\begin{equation*}
\begin{array}{cccc}
\mathcal{D}_{p}\subseteq & \mathcal{D} & \subseteq & \mathcal{A} \\ 
& \cap &  & \cap \\ 
& \mathcal{B} & \subseteq & \tprod\limits_{k}\mathcal{M}_{n_{k}}\left( 
\mathbb{C}\right) /\sum_{k}\mathcal{M}_{n_{k}}\left( \mathbb{C}\right)%
\end{array}%
\end{equation*}%
commutes. By Theorem \ref{16}, we obtain that $\mathcal{A}\ast _{\mathcal{D}%
_{p}}\mathcal{B}$ is an MF algebra for each $p\geq 1$. Now it follows from
Theorem \ref{Theorem 3.1} that $\mathcal{A}\ast _{\mathcal{D}}\mathcal{B}$
is an MF algebra.
\end{proof}

Since every AF algebra has a faithful tracial state, we are able to consider
the case when $\mathcal{A}$, $\mathcal{B}$ and $\mathcal{D}$ are all AF
algebras and give a sufficient condition in terms of faithful tracial states.

\begin{theorem}
{\label{Theorem 3.3}} Suppose that $\mathcal{A}\supset \mathcal{D}\subset 
\mathcal{B}$ are unital inclusions of AF C$^{\ast }$-algebras. If there are
faithful tracial states $\tau _{\mathcal{A}}$ and $\tau _{\mathcal{B}}$ on $%
\mathcal{A}$ and $\mathcal{B}$ respectively, such that 
\begin{equation*}
\tau _{\mathcal{A}}(x)=\tau _{\mathcal{B}}(x),\qquad \forall \ x\in \mathcal{%
D},
\end{equation*}%
then $\mathcal{A}\ast _{\mathcal{D}}\mathcal{B}$ is an MF algebra.
\end{theorem}

\begin{proof}
Assume that $\{x_{n}\}_{n=1}^{\infty }\subseteq \mathcal{A}$, $%
\{y_{n}\}_{n=1}^{\infty }\subseteq \mathcal{B}$ and $\{z_{n}\}_{n=1}^{\infty
}\subseteq \mathcal{D}$ are families of generators in $\mathcal{A}$, $%
\mathcal{B}$ and $\mathcal{D}$ respectively. Note that $\mathcal{A},\mathcal{%
B}$ and $\mathcal{D}$ are AF algebras. For each $N\in \mathbb{N}$, there are
finite dimensional C$^{\ast }$-subalgebras $\mathcal{D}_{N}\subseteq 
\mathcal{D}$, $\mathcal{A}_{N}\subseteq \mathcal{A}$ and $\mathcal{B}%
_{N}\subseteq \mathcal{A}$ such that 
\begin{equation}
\max_{1\geq n\leq N}\{dist(x_{n},\mathcal{A}_{N}),dist(y_{n},\mathcal{B}%
_{N}),dist(z_{n},\mathcal{D}_{N})\}\leq \frac{1}{N}  \tag{5.43}
\label{equation 3.4}
\end{equation}%
and 
\begin{equation}
\mathcal{A}_{N}\supset \mathcal{D}_{N}\subset \mathcal{B}_{N}  \tag{5.44}
\label{equation 3.5}
\end{equation}%
Note that $\tau _{\mathcal{A}}(x)=\tau _{\mathcal{B}}(x),\ \forall \ x\in 
\mathcal{D}.$ From the argument in the proof of Theorem 4.2 \cite{[ADRL]},
there are rational faithful tracial states on $\mathcal{A}_{N}$ and $%
\mathcal{B}_{N}$ such that their restrictions on $\mathcal{D}_{N}$ agree$.$
This implies that there is a positive integer $k_{N}$ such that 
\begin{equation}
\mathcal{M}_{k_{N}}(\mathbb{C)\supseteq \mathcal{A}}_{N}\mathbb{\supseteq 
\mathcal{D}}_{N}\mathbb{\subseteq \mathcal{B}}_{N}\mathbb{\subseteq \mathcal{%
M}}_{k_{N}}\mathbb{(C)}.  \tag{5.45}  \label{equation 3.6}
\end{equation}%
Combining \ref{equation 3.4}, \ref{equation 3.5} and \ref{equation 3.6}, we
know that there is a sequence of positive integers $\{k_{N}\}_{N=1}^{\infty
} $ such that the diagram 
\begin{equation*}
\begin{array}{ccc}
\mathcal{D} & \subseteq & \mathcal{A} \\ 
\cap &  & \cap \\ 
\mathcal{B} & \subseteq & \tprod\limits_{N}\mathcal{M}_{k_{N}}\left( \mathbb{%
C}\right) /\sum_{N}\mathcal{M}_{k_{N}}\left( \mathbb{C}\right)%
\end{array}%
\end{equation*}%
commutes. By Theorem \ref{Theorem 3.2}, we obtain that $\mathcal{A}\ast _{%
\mathcal{D}}\mathcal{B}$ is an MF algebra.
\end{proof}

It is well-known that the tracial state on each UHF algebra is unique.
Therefore we can restate Theorem \ref{Theorem 3.3} when $\mathcal{A}$, $%
\mathcal{B}$ are both UHF algebras and $\mathcal{D}$ is an AF algebra. A
necessary and sufficient condition can be given in this case.

\begin{corollary}
Suppose that $\mathcal{A\supseteq D\subseteq B}$ are unital inclusions of
C*-algebras where $\mathcal{A}$, $\mathcal{B}$ are UHF algebras and $%
\mathcal{D}$ is an AF algebra. Then $\mathcal{A}\underset{\mathcal{D}}{%
\mathcal{\ast }}\mathcal{B}$ is an MF algebra if and only if 
\begin{equation*}
\tau _{\mathcal{A}}\left( z\right) =\tau _{\mathcal{B}}\left( x\right) \ 
\text{for\ each\ }z\in \mathcal{D},
\end{equation*}%
where $\tau _{\mathcal{A}}$ and $\tau _{\mathcal{B}}$ are faithful tracial
states on $\mathcal{A}$ and $\mathcal{B}$ respectively.
\end{corollary}

\begin{proof}
From the fact that every MF algebra has a tracial state and the tracial
state on UHF algebra is unique and faithful, one direction of the proof is
obvious. Another direction is followed by applying Theorem \ref{Theorem 3.3}.
\end{proof}

{\small \bigskip \allowbreak }\newpage

Junhao Shen

Department of Mathematics and Statistics

Univerisity of New Hampshire

Durham, NH, 03824

Email: jog2@cisunix.unh.edu\bigskip 
\begin{equation*}
\end{equation*}

Qihui Li

Department of Mathematics and Statistics

Univerisity of New Hampshire

Durham, NH, 03824

Email: qme2@cisunix.unh.edu\bigskip

and

Department of Mathematics

East China University of Science and Technology

Shanghai, China Meilong Road 130, 200237

\end{document}